\documentclass[11pt]{amsart}
\usepackage[a4paper,margin=26mm]{geometry}
\usepackage{amsmath,amssymb,amsthm,booktabs,array,enumitem}
\usepackage[hidelinks]{hyperref} % 修改处：去掉了 unicode 选项
\allowdisplaybreaks[2]
\newtheorem{theorem}{Theorem}[section]
\newtheorem{lemma}[theorem]{Lemma}
\newtheorem{proposition}[theorem]{Proposition}
\newtheorem{corollary}[theorem]{Corollary}

\numberwithin{equation}{section}
\DeclareMathOperator{\supp}{supp}
\newcommand{\height}{\mathsf h}
\newcommand{\vnum}{\mathsf v}
\newcommand{\caseheading}[2]{%
	\par\addvspace{0.5\baselineskip}%
	\noindent\textit{Case~#1. #2.}\enspace\ignorespaces
}
\title{The zero-sum constant \(\mathsf{s}_p^*\) and inverse problem}

\author{Cong Lin} 
\date{}
\keywords{Zero-sum sequences, inverse zero-sum problems, cyclic groups, restricted sumsets}
\begin{document}
	\maketitle

\begin{abstract}
Let $p$ be an odd prime. We denote by $\mathsf{s}_{p}^{*}$ the smallest integer $l$ such that,
out of any given $l$ integers coprime with $p$,
one can select $p$ integers satisfying that
the sum of the $p$ integers is a multiple of $p$,
but not a multiple of $p^{2}$.
It was conjectured that $\mathsf{s}_p^*=2p+1$. In this paper we prove this conjecture
and solve the corresponding inverse problem.
More precisely, we characterize all sequences of $2p$ integers coprime to $p$
in which every $p$-term subsequence with sum a multiple of $p$
also has sum a multiple of $p^2$.

\end{abstract}

\section{Introduction}

Two classical constants distinguish between finding a nonempty
zero-sum subsequence and finding one of prescribed length. For a finite
abelian group $H$, the Davenport constant $\mathsf D(H)$ is the least
integer $l$ such that every sequence of length $l$ over $H$ has a
nonempty zero-sum subsequence. The Erd\H{o}s--Ginzburg--Ziv constant
$\mathsf s(H)$ is defined in the same way, with the additional requirement
that the subsequence have length $\exp(H)$, the exponent of $H$.
Let $C_n$ denote the cyclic group of order $n$. For a prime $p$,
we have $\mathsf D(C_p)=p$, and the
Erd\H{o}s--Ginzburg--Ziv theorem gives $\mathsf s(C_p)=2p-1$.
Thus among any $2p-1$ integers one can find $p$ terms whose sum is
divisible by $p$.

In their study of zero-sum subsequences of cross number $1$,
Gao, Hui, Jiang, Li and Wang asked for a $p$-term sum divisible
by $p$ but not by $p^2$ when the terms are coprime to $p$
\cite[Conjecture~3.2]{GHJLW}. Gao, Jiang, Lei, Lin and Yang formulated
this question in terms of the following constant.

\par\medskip
\noindent\textbf{Definition.}
For a prime $p$, let $\mathsf{s}_p^*$ be the smallest integer $l$
such that every sequence of $l$ integers coprime to $p$ has a
subsequence of length $p$ whose sum is divisible by $p$ but not by
$p^2$ \cite[Definition~1.1]{GJLLY}.

\medskip
For odd $p$, the sequence $1^{p-1}(1-p)(-1)^{p-1}(p-1)$ in
\cite{GHJLW} gives $\mathsf{s}_p^*\ge2p+1$, and the conjecture asserts
that equality holds. The upper bound $\mathsf{s}_p^*\le4p-3$ follows
from Reiher's theorem $\mathsf{s}(C_p\oplus C_p)=4p-3$
\cite{Reiher}. Represent each term modulo $p^2$ as $u+pv$ with
$1\le u\le p-1$ and $0\le v\le p-1$, and apply the theorem to
the pairs $(u,v)$ modulo $p$ (see \cite{GJLLY}). In the selected
$p$ terms, the sum of the $u$-coordinates is a positive multiple
of $p$ smaller than $p^2$, and the sum of the $v$-coordinates
is a multiple of $p$.
Gao, Jiang, Lei, Lin and Yang improved the upper bound to $3p-2$
\cite[Theorem~1.3]{GJLLY} and reported that equality with $2p+1$
had been verified for $p=3,5,7$.

Related problems differ in the permitted lengths or in the modulus.
In a preprint, Gao, Jiang and Mu showed that for a prime power $n$,
every sequence of $2n-1$ integers coprime to $n$ has a nonempty
subsequence of length at most $n$ with sum $u$ satisfying
$\gcd(u,n^2)=n$. They also determined the corresponding extremal
sequences \cite{GJM}. For $n=p$ the condition on the sum is the same
as in the definition of $\mathsf{s}_p^*$, with the length allowed
to be smaller than $p$. Sun considers versions with exactly $n$
terms and general integer modulus $n$ in the preprint \cite{SunLocal}.

The defining condition for $\mathsf{s}_p^*$ depends only on the
residues modulo $p^2$. Let $G=\mathbb Z/p^2\mathbb Z$, whose additive
group is $C_{p^2}$. The ideal $pG=\{pg:g\in G\}$ is an additive
subgroup of order $p$, and $G^\times=G\setminus pG$ is the set of
units of $G$. We identify $C_p$ with the additive group of
$\mathbb Z/p\mathbb Z$. The canonical homomorphism
$\pi\colon G\to C_p$ is reduction modulo $p$ and has kernel $pG$.
For a sequence $T$ over $G$, let
$\overline T$ be the sequence obtained by applying $\pi$ to each term of $T$.

For a sequence $T$, let $|T|$ be its length and $\sigma(T)$ its sum.
A sequence $S$ of length $2p$ over $G^\times$ is called exceptional
if it has no subsequence $T$ of length $p$ with
$\sigma(T)\in pG\setminus\{0\}$.

The following theorem classifies the exceptional sequences.
Its corollary proves $\mathsf{s}_p^*=2p+1$ for all odd primes, so the
classification gives all sequences of length $\mathsf{s}_p^*-1$
that fail the defining condition.

\begin{theorem}\label{thm:main}
Let $p$ be an odd prime and let $S$ be a sequence of length $2p$ over $G^\times$.
Then $S$ is exceptional if and only if it has one of the following forms.
\begin{enumerate}[label=\textup{(\roman*)},ref=\roman*]
\item\label{item:two-classes}
$S=AB$, where $r,t$ are distinct nonzero elements of $C_p$ and
$$
 |A|=|B|=p,\qquad \overline A=r^p,\qquad \overline B=t^p,\qquad
 \sigma(A)=\sigma(B)=0.
$$
\item\label{item:four-terms}
$S=a^{p-1}b^{p-1}xy$, where $a,b,x,y\in G^\times$, $a-b\in G^\times$,
and there exists an integer $i$ with $1\le i\le p-2$ such that
$$
 x=-ia-(p-1-i)b,\qquad y=-(p-1-i)a-ib.
$$
The terms $x$ and $y$ may coincide.
\end{enumerate}
\end{theorem}

\begin{corollary}\label{cor:direct}
For every odd prime $p$, one has $\mathsf{s}_p^*=2p+1$.
\end{corollary}

The proof begins with the case of two distinct terms with the same
image in $C_p$. A common completion to a $p$-term subsequence gives
a contradiction unless the sequence is concentrated in two nonzero
residue classes, each containing $p$ terms. When $\pi$ is
injective on the support, we first treat sequences with large
multiplicities. In the remaining cases, we distribute the copies of
each value among layers and use polynomial interpolation modulo $p$
and $p^2$ to bound the corresponding sumsets. These bounds allow us
to lift selected sumsets to integer intervals and then complete their
sums to $p$ terms.

\section{Notation and preliminaries}\label{sec:preliminaries}

\subsection{Sequences and their images}

We denote by $\mathbb N$ the set of positive integers and by
$\mathbb N_0=\mathbb N\cup\{0\}$ the set of nonnegative integers.
For real numbers $a,b$, we set $[a,b]=\{n\in\mathbb Z:a\le n\le b\}$.

Let $H$ be a finite additive abelian group. We denote by
$\mathcal F(H)$ the free abelian monoid with basis $H$. Its elements
are called sequences over $H$. We identify $g\in H$ with the
one-term sequence $g$. We denote concatenation by juxtaposition,
and $g^j$ denotes the sequence of $j$ copies of $g$. Unless
otherwise stated, $\prod$ denotes concatenation of an indexed
family of sequences.
For $n\in\mathbb Z$ and $g\in H$, $ng$ denotes the integer multiple
of $g$ in $H$. For $S\in\mathcal F(H)$ and $x\in H$,
let $\vnum_x(S)\in\mathbb N_0$ be the multiplicity of $x$ in $S$.
We write $T\mid S$ if $\vnum_x(T)\le\vnum_x(S)$ for all $x\in H$.
Then $T$ is a subsequence of $S$, and $ST^{-1}$ denotes the sequence
obtained by deleting $T$ from $S$. The identity of $\mathcal F(H)$
is the empty sequence.

For a sequence
$$
 S=g_1\cdots g_l
  =\prod_{x\in H}x^{\vnum_x(S)}\in\mathcal F(H),
$$
we use the following notation.
\begin{itemize}[leftmargin=*,nosep]
\item $|S|=l=\sum_{x\in H}\vnum_x(S)$ is the length of $S$.
\item $\sigma(S)=\sum_{i=1}^{l}g_i=\sum_{x\in H}\vnum_x(S)x$
      is the sum of $S$.
\item $\supp(S)=\{x\in H:\vnum_x(S)>0\}$ is the support of $S$.
\item $\height(S)=\max\{\vnum_x(S):x\in H\}$ is the maximum multiplicity
      in $S$.
\item $S$ is a zero-sum sequence if $\sigma(S)=0$. A nonempty zero-sum
      sequence is minimal if it has no nonempty proper zero-sum subsequence.
      A sequence is zero-sum free if it has no nonempty zero-sum subsequence.
\item $\Sigma_j(S)=\{\sigma(T):T\mid S,\ |T|=j\}$ for $j\in[0,|S|]$.
      We denote by $\Sigma(S)$ the set of all nonempty subsequence sums.
\end{itemize}

For $g\in G$, let $\overline g=\pi(g)$. For a sequence
$T=t_1\cdots t_k$ over $G$, let
$\overline T=\overline t_1\cdots\overline t_k$.
Since $\pi$ is a homomorphism with kernel $pG$, we have
$\sigma(\overline T)=\pi(\sigma(T))$ and
$\sigma(\overline T)=0$ if and only if $\sigma(T)\in pG$.
For a subset $A\subseteq G$, let $\pi(A)=\{\pi(a):a\in A\}$.

Equivalently, an exceptional sequence $S$ satisfies
\begin{equation}\label{eq:exceptional}
 |T|=p,\quad T\mid S,\quad \sigma(\overline T)=0
 \quad\Longrightarrow\quad\sigma(T)=0.
\end{equation}

Since $\mathsf{s}(C_p)=2p-1$, $\overline S$ has a $p$-term zero-sum subsequence for every
exceptional sequence $S$.

\subsection{Sumset bounds and long zero-sum sequences}

For subsets $A,B$ of an additive abelian group, write
$$
 -A=\{-a:a\in A\},\qquad A+B=\{a+b:a\in A,\ b\in B\}.
$$
The sumset of $A$ and $B$ is $A+B$, and their difference set is
$A-B=A+(-B)$. For a group element $x$, we denote $\{x\}+A$ by
$x+A$. For a positive integer $r$, the sumset of $r$ copies of
$A$ is denoted by $rA$. For an integer set $I$ and a group
element $g$, we use $Ig=\{ig:i\in I\}$.

We use the following form of the Cauchy--Davenport theorem, obtained
by iterating the two-set inequality.

\begin{lemma}[{\cite{Cauchy,Davenport}}]\label{lem:cauchy-davenport}
Let $p$ be a prime and let $A_1,\ldots,A_r$ be nonempty subsets of
$C_p$, where $r\ge1$. Then
\begin{equation}\label{eq:cd}
 |A_1+\cdots+A_r|\ge\min\left\{p,1+\sum_{j=1}^{r}(|A_j|-1)\right\}.
\end{equation}
\end{lemma}

We also use the following form of Vosper's theorem~\cite{Vosper}.
If $A,B\subseteq C_p$, $|A|,|B|\ge2$, and
\begin{equation}\label{eq:vosper}
 |A+B|=|A|+|B|-1\le p-2,
\end{equation}
then $A$ and $B$ are arithmetic progressions with the same common
difference.

We use Savchev and Chen's theorem on long zero-sum free sequences
and its consequence for long minimal zero-sum sequences.

\begin{lemma}[{\cite[Theorem~8 and Section~5]{SC}}]\label{lem:long-normal-form}
Let $V$ be a sequence over $C_p$.
\begin{enumerate}[label=\textup{(\roman*)},nosep]
\item If $V$ is zero-sum free and $|V|>p/2$, then for some
generator $g$ of $C_p$,
$$
 V=(a_1g)\cdots(a_{|V|}g),\qquad
 a_i\in\mathbb N,\qquad \sum_{i=1}^{|V|}a_i\le p-1.
$$
\item If $V$ is minimal zero-sum and $|V|>(p+2)/2$, then
$V$ has such a representation with $\sum_{i=1}^{|V|}a_i=p$.
\end{enumerate}
\end{lemma}

\begin{lemma}[{\cite[Section~6]{SC}}]\label{lem:coefficient-multiplicity}
Let $V=(a_1g)\cdots(a_l g)$, where $g$ is a generator
of $C_p$ and $a_1,\ldots,a_l\in\mathbb N$. If
$N=\sum_{i=1}^{l}a_i$, then $\vnum_g(V)\ge2|V|-N$.
\end{lemma}

For a minimal zero-sum sequence $U$ over $C_p$, we have
\begin{equation}\label{eq:minsum-high}
 |\{0\}\cup\Sigma(U)|\ge\min\{p,2|U|-1\}.
\end{equation}
Write $U=u_1\cdots u_k$, set $s_0=0$, and put
$s_i=\sum_{l=1}^i u_l$ for $1\le i<k$. Let
$P=\{s_0,\ldots,s_{k-1}\}$. For $0\le i<j<k$, the subsequence
$W=u_{i+1}\cdots u_j$ is nonempty and proper. We have
$$
 s_j-s_i=\sigma(W)\ne0,\qquad
 s_i-s_j=\sigma(UW^{-1}),
$$
where the first assertion follows from minimality and the second from
$\sigma(U)=0$. Thus $|P|=k$ and $P-P\subseteq\{0\}\cup\Sigma(U)$.
Applying Lemma~\ref{lem:cauchy-davenport} to $P$ and $-P$ gives
\eqref{eq:minsum-high}.

\begin{lemma}\label{lem:fixed-length-gaps}
Let $z_1\le\cdots\le z_n$ be integers with
$z_{i+1}-z_i\le d$ for $1\le i<n$, where $d\ge1$ is an integer.
For $k\in[0,n]$, let $Y$ be the set of sums of $k$ distinct positions.
Successive elements of $Y$ differ by at most $d$.
Moreover, if $P$ is a nonempty integer interval and $Z$ is a
nonempty finite set of integers whose successive elements differ
by at most $|P|$, then $P+Z$ is an integer interval.
\end{lemma}

\begin{proof}
For $k=0$ or $k=n$, the set $Y$ is a singleton. For $1\le k<n$,
start with the positions $1,\ldots,k$. Move the selected positions
from right to left to their destinations $n-k+1,\ldots,n$, advancing
one position at a time. This gives a nondecreasing chain of sums
from $\min Y$ to $\max Y$ with increments at most $d$. Hence each
gap in $Y$ is at most $d$. For the second assertion, consecutive
translates $P+z$ for $z\in Z$ overlap or are adjacent, so their
union $P+Z$ is an interval.
\end{proof}

We use the following special case of Grynkiewicz's classification
of sequences with a unique zero-sum length.

\begin{lemma}[{\cite[Theorem~3.4]{Gry}}]\label{lem:unique-zero-sum-length}
Let $p$ be an odd prime and let $T$ be a sequence over $C_p$ with
$|T|=p$ and $\height(T)\le p-2$. If there is a unique $r\in[1,p]$
with $0\in\Sigma_r(T)$, then there is a generator $g$ of $C_p$ such that
$$
 r=(p+1)/2,\qquad T=g^{p-2}(rg)^2.
$$
\end{lemma}

\begin{lemma}\label{lem:lift-coefficient}
Let $A=UV$ be a sequence over $G$ with $\sigma(\overline V)=0$.
Suppose that all nonempty proper subsequences of $A$ with sum
in $pG$ have the same sum in $G$. Let $j\ge1$ and $w,q\in G$
with $wq^j\mid V$ and $\overline w=j\overline q$.
Set $V_0=V(wq^j)^{-1}$. If
$$
 \bigl(j\overline q+(\{0\}\cup\Sigma(\overline{V_0}))\bigr)
 \cap \bigl(\Sigma(\overline U)\setminus\{0\}\bigr)\ne\varnothing,
$$
then $w=jq$.
\end{lemma}

\begin{proof}
Choose a nonempty subsequence $D\mid U$ and a subsequence
$W\mid V_0$, possibly empty, such that
$\sigma(\overline D)=j\overline q+\sigma(\overline W)\ne0$.
Consider
$$
 T_1=DV(wW)^{-1},\qquad T_2=DV(q^jW)^{-1}.
$$
Both retain $D$ and omit terms of $V$, so they are nonempty
proper subsequences of $A$. By $\sigma(\overline V)=0$,
$$
 \sigma(\overline{T_1})=\sigma(\overline{T_2})
 =\sigma(\overline D)-j\overline q-\sigma(\overline W)=0.
$$
The hypothesis on $A$ now gives $0=\sigma(T_1)-\sigma(T_2)=jq-w$.
\end{proof}

\section{Noninjectivity and large multiplicities}\label{sec:high}

Unless otherwise stated, from this section through the proof of
Theorem~\ref{thm:main}, $S$ denotes a sequence over $G^\times$ of
length $2p$, where $p$ is an odd prime.

If $\height(S)\ge p$, then some unit $a$ occurs at least $p$ times in $S$.
The $p$-term subsequence $a^p$ has sum $pa\in pG\setminus\{0\}$. Henceforth assume $\height(S)\le p-1$.

Gao et al.~\cite[Lemma~3.4]{GHJLW} showed that if a sequence of
length $2p+1$ over $G^\times$ has no $p$-term subsequence with
sum in $pG\setminus\{0\}$, then $\pi$ is injective on its support.
The following proposition classifies the exceptional sequences of
length $2p$ for which $\pi$ fails to be injective on the support.

\begin{proposition}\label{prop:collision}
If $\pi$ is not injective on $\supp(S)$, then $S$ is exceptional if and only if it has the form in Theorem~\ref{thm:main}\textup{(\ref{item:two-classes})}.
\end{proposition}
\begin{proof}
Assume that $S$ is exceptional. At most $p$ terms of $S$ have any
given image under $\pi$. Indeed, if $L\mid S$ consists of $p+1$ terms
with the same image, then $\sigma(\overline{Lx^{-1}})=p\overline x=0$
for every term $x$ of $L$.
Thus $\sigma(L)-x=0$ by \eqref{eq:exceptional}, so all terms of $L$
are equal, contradicting $\height(S)\le p-1$.

Choose distinct terms $a,b$ with $\overline a=\overline b=r$ and let
$S_0=S(ab)^{-1}$. Suppose $\height(\overline{S_0})\le p-1$. List the
terms of $\overline{S_0}$ as $r_1,\ldots,r_{2p-2}$ with equal values
consecutive, and set $P_i=\{r_i,r_{i+p-1}\}$ for $1\le i\le p-1$.
Each $P_i$ has two elements, since $r_i=r_{i+p-1}$ would force
$p$ consecutive terms to be equal. By Lemma~\ref{lem:cauchy-davenport},
$$
 p\le |P_1+\cdots+P_{p-1}|
   \le|\Sigma_{p-1}(\overline{S_0})|\le p.
$$
Thus $\Sigma_{p-1}(\overline{S_0})=C_p$. Choose $W\mid S_0$ with
$|W|=p-1$ and $\sigma(\overline W)=-r$.
By \eqref{eq:exceptional}, $\sigma(aW)=\sigma(bW)=0$, and hence
$a=b$, a contradiction.

It follows that $S_0$ has a subsequence $A$ of length $p$ whose
terms have a common image $t\ne r$, since at most $p$ terms of $S$
have image $r$ and two were deleted. The sum of $A$ is zero, and
$\height(S)\le p-1$ ensures that $A$ contains two distinct terms.
Applying the preceding argument to these terms gives a second
subsequence $B$ of length $p$ and sum zero, whose terms have a
common image different from $t$. Thus $S=AB$, and the terms of $B$
have image $r$ since $B$ contains $a$ and $b$.

Conversely, if $T\mid AB$ uses $j$ terms of $A$ and $p-j$ terms of
$B$, then
$$
 \sigma(\overline T)=jt+(p-j)r=j(t-r).
$$
As $0\le j\le p$ and $t\ne r$, this vanishes only for $j=0$ or $j=p$.
In these cases $T$ is $B$ or $A$, respectively, and $\sigma(T)=0$.
\end{proof}

It remains to classify exceptional sequences satisfying
\begin{equation}\label{eq:common}
\pi\text{ is injective on }\supp(S).
\end{equation}
Choose $a\in\supp(S)$ with $\vnum_a(S)=\height(S)$, set
$m=p-\height(S)$, and let $R=Sa^{-(p-m)}=x_1\cdots x_{p+m}$.
Thus $R$ is obtained by deleting all occurrences of $a$ from $S$. Set
$$
 Q=\prod_{j=1}^{p+m}(x_j-a),\qquad c=-pa.
$$
Since $\supp(Q)=\supp(R)-\{a\}$, \eqref{eq:common} gives
injectivity on $\supp(Q)$ and $\vnum_0(\overline Q)=0$.
In particular, $\height(\overline Q)=\height(Q)\le p-m$. For
$I\subseteq[1,p+m]$ with $m\le |I|\le p$, let
$$
 B=\prod_{i\in I}(x_i-a),\qquad T=a^{p-|I|}\prod_{i\in I}x_i.
$$
Since $p-|I|\le p-m$, we have $T\mid S$, $|T|=p$, and
$\sigma(T)=pa+\sigma(B)$. If $S$ is exceptional, then
\eqref{eq:exceptional} gives
\begin{equation}\label{eq:fixed}
	B\mid Q,\quad m\le |B|\le p,\quad \sigma(B)\in pG
	\quad\Longrightarrow\quad \sigma(B)=c.
\end{equation}

\subsection{Maximum multiplicity \texorpdfstring{$p-1$}{p-1}}

We first exclude exceptional sequences with exactly one value of
multiplicity $p-1$ and then characterize those with two such values.

\begin{proposition}\label{prop:base}
Under the injectivity assumption~\eqref{eq:common}, if exactly one value has multiplicity $p-1$ in $S$, then $S$ has a subsequence $T$ with $|T|=p$ and $\sigma(T)\in pG\setminus\{0\}$.
\end{proposition}
\begin{proof}
Assume that $S$ is exceptional. Then $|Q|=p+1$ and
$\height(Q)\le p-2$. Define
$$
 A=
 \begin{cases}
 Q, & \sigma(\overline Q)=0,\\
 Qz,\quad z=2c-\sigma(Q), & \sigma(\overline Q)\ne0.
 \end{cases}
$$
In the second case, $\overline z=-\sigma(\overline Q)\ne0$. Thus in both
cases $\sigma(\overline A)=0$ and $\vnum_0(\overline A)=0$.
We claim that
\begin{equation}\label{eq:fixed-A}
 B\mid A,\quad 1\le |B|<|A|,\quad \sigma(B)\in pG
 \quad\Longrightarrow\quad \sigma(B)=c.
\end{equation}
If $A=Q$, the claim is \eqref{eq:fixed} with $m=1$.

Suppose $A=Qz$, and let $B$ satisfy the conditions of \eqref{eq:fixed-A}. Then
$\sigma(AB^{-1})=2c-\sigma(B)\in pG$. Neither $B$ nor $AB^{-1}$
can equal $Q$, since $\sigma(Q)\notin pG$. Exactly one of the two
subsequences omits the position of $z$. By \eqref{eq:fixed} that
subsequence has sum $c$. Since $\sigma(B)+\sigma(AB^{-1})=2c$,
both sums are $c$.

We next show that $\pi$ is injective on $\supp(A)$. It suffices to check
the case with an added term. If $x,y$ are two terms of $A$ with
$\overline x=\overline y$, put $A'=A(xy)^{-1}$. Then $|A'|=p$,
and applying Lemma~\ref{lem:cauchy-davenport} to $\{0,\overline v\}$ for the positions
$v$ of $A'$ gives
$$
 \{0\}\cup\Sigma(\overline{A'})=C_p.
$$
Choose $W\mid A'$ with $\sigma(\overline W)=-\overline x$.
We have
$$
 1\le |xW|=|yW|\le p+1<|A|,\qquad
 \sigma(\overline{xW})=\sigma(\overline{yW})=0.
$$
By \eqref{eq:fixed-A}, $\sigma(xW)=\sigma(yW)=c$, so $x=y$.

Take a nonempty subsequence $U\mid A$ of minimum length subject to
$\sigma(\overline U)=0$. Since $\mathsf D(C_p)=p<|A|$, we have
$|U|\le p$. Let $V=AU^{-1}$, $r=|U|$, and $\varepsilon=|A|-p$.
By the choice of $U$, $\overline U$ is minimal zero-sum.
Also, $V$ is a nonempty proper subsequence of $A$ with
$\sigma(\overline V)=\sigma(\overline A)-\sigma(\overline U)=0$.
By \eqref{eq:fixed-A}, $\sigma(U)=\sigma(V)=c$.

Suppose that $\overline V$ is not minimal zero-sum. Then there is
a subsequence $W\mid V$ with $1\le |W|<|V|$ and
$\sigma(\overline W)=0$. Since $\sigma(\overline V)=0$, we also
have $\sigma(\overline{VW^{-1}})=0$. Apply \eqref{eq:fixed-A}
to $W$ and $VW^{-1}$ to obtain
$$
 c=\sigma(V)=\sigma(W)+\sigma(VW^{-1})=2c,
$$
contrary to $c\ne0$. Hence $\overline V$ is minimal zero-sum.

Since every term of $\overline A$ is nonzero and $U$ has minimum
length, we have $2\le r\le |V|$. Thus $2r\le |A|\le p+2$.
As $p$ is odd, it follows that $r\le(p+1)/2$.

We distinguish the cases according to whether $|V|>(p+2)/2$.

\caseheading{1}{$|V|\le(p+2)/2$}
Since $|V|=p+\varepsilon-r$ and $r\le(p+1)/2$, we have
$\varepsilon=1$ and $r=(p+1)/2$. Thus $A=Q$.
Each nonempty proper zero-sum subsequence of $\overline Q$
and its complement have length at least $r$, so both have length $r$.
Delete any term $y$ of $Q$. Since $|\overline{Qy^{-1}}|=p=\mathsf D(C_p)$,
the sequence $\overline{Qy^{-1}}$ has a nonempty zero-sum subsequence,
and $r$ is the only possible length. Lemma~\ref{lem:unique-zero-sum-length}
gives $\overline{Qy^{-1}}=g^{p-2}(rg)^2$ for some generator $g$.
Therefore $\sigma(\overline{Qy^{-1}})=-g$ and $\overline y=g$.
This gives $\vnum_g(\overline Q)\ge p-1$, a contradiction.

\caseheading{2}{$|V|>(p+2)/2$}
By Lemma~\ref{lem:long-normal-form}\textup{(ii)}, choose a generator $g$
of $C_p$ so that the terms of $\overline V$ have the form $jg$ with
$j\in[1,p-1]$ and these integer coefficients sum to $p$.
By Lemma~\ref{lem:coefficient-multiplicity} and injectivity on $\supp(A)$,
the term $q$ with $\overline q=g$ satisfies $\vnum_q(V)\ge2|V|-p\ge3$,
so $q\mid Q$. In particular, $a+q\in\supp(S)\subseteq G^\times$.

We first prove that a term $w\mid V$ of coefficient $j\in[2,r-1]$
satisfies $w=jq$. Applying Lemma~\ref{lem:coefficient-multiplicity} to
$\overline{Vw^{-1}}$ gives
$$
 \vnum_q(V)\ge2(|V|-1)-(p-j)
   =p+2\varepsilon-2r+j-2\ge j+1.
$$
Let $V_0=V(wq^j)^{-1}$, and set
$$
 X=jg+(\{0\}\cup\Sigma(\overline{V_0})),\qquad
 Y=\{0\}\cup\Sigma(\overline U).
$$
The positive coefficients of $\overline{V_0}$ sum to $p-2j$, so
$X\subseteq[j,p-j]g$ and $0\notin X$. The set
$\{0\}\cup\Sigma(\overline{V_0})$ is the sum of the sets
$\{0,\overline v\}$, one for each term $v$ of $V_0$. Each has
two elements. By Lemma~\ref{lem:cauchy-davenport},
$|X|\ge |V_0|+1=p+\varepsilon-r-j$.
The same bound holds for empty $V_0$, when $X=\{jg\}$. Since
$2r-1\le p$, \eqref{eq:minsum-high} gives $|Y|\ge2r-1$. Hence
$$
 |X|+|Y|-1\ge p+\varepsilon+r-j-2\ge p.
$$
Applying Lemma~\ref{lem:cauchy-davenport} to $X$ and $-Y$ gives $X-Y=C_p$.
Thus $0\in X-Y$, so $X\cap Y\ne\varnothing$. Since $0\notin X$,
\eqref{eq:fixed-A} and Lemma~\ref{lem:lift-coefficient} give $w=jq$.

\caseheading{2.1}{Every coefficient of $\overline V$ is less than $r$}
All terms of $V$ equal their corresponding multiples of $q$, so
$pq=\sigma(V)=c$, contradicting that $a+q$ is a unit.

\caseheading{2.2}{Some coefficient of $\overline V$ is at least $r$}
A term $w$ of coefficient $j\ge r$ satisfies
$$
 r-1\le j-1\le p-|V|=r-\varepsilon.
$$
Thus $\varepsilon=1$, $j=r$, and $V=q^{p-r}w$. Since
$|V|>(p+2)/2$, we have $r\le(p-1)/2$. For $V_0=q^{p-2r}$, the
translated sumset in Lemma~\ref{lem:lift-coefficient} is $[r,p-r]g$.
If it met $\{0\}\cup\Sigma(\overline U)$, \eqref{eq:fixed-A} and
the lemma would give $w=rq$ and hence $pq=\sigma(V)=c$, contrary
to $a+q$ being a unit. Thus \eqref{eq:minsum-high} gives
$$
 \{0\}\cup\Sigma(\overline U)=[-(r-1),r-1]g.
$$

Each term of $\overline U$ has the form $jg$ for a unique
$j\in[-(r-1),r-1]\setminus\{0\}$. A term $u\mid U$ with
$\overline u=-jg$ and $1\le j\le r-2$
would give $uq^j\mid A$ with $\sigma(\overline{uq^j})=0$ and
$|uq^j|<r$. Hence the only possible negative coefficient in $U$
is $-(r-1)$. If $z_0$ is a term with this coefficient,
\eqref{eq:fixed-A} applied to $z_0q^{r-1}$ gives
$z_0=c-(r-1)q$. If $u\mid U$ has $\overline u=jg$ with $1\le j\le r-1$,
then $Uu^{-1}q^j\mid A$ has image sum zero and length at most
$2r-2<|A|$. By \eqref{eq:fixed-A},
$$
 \sigma(Uu^{-1}q^j)=\sigma(U)=c,\qquad u=jq.
$$
Let $l$ be the number of negative coefficients and $P$ the sum of
the positive coefficients.

\caseheading{2.2.1}{$l=0$}
Choose a subsequence $D\mid U$ of minimum length whose
coefficient sum $d$ is at least $r$. Removing any term lowers the
sum below $r$. Since the positive coefficients are at most $r-1$,
$$
 r\le d\le 2r-2<p,\qquad
 dg\in [-(r-1),r-1]g.
$$
It follows that $d\ge p-r+1$. Since $0<d<p$ and
$\sigma(\overline U)=0$, the subsequence $D$ is proper in $U$.
Let $T=Dq^{p-d}$. Since $p-d\le p-r$, we have $T\mid A$.
Moreover, $1\le |T|\le 2r-2<|A|$ and
$$
 \sigma(T)=dq+(p-d)q=pq\in pG.
$$
By \eqref{eq:fixed-A}, $\sigma(T)=c$, so $pq=c$, a contradiction.

\caseheading{2.2.2}{$l=1$}
Here $\sigma(U)=Pq+z_0=c$, and hence
$$
 (P-r+1)q=0,\qquad 0\le P-r+1<p^2.
$$
Since $q$ has order $p^2$, we have $P=r-1$. The $r-1$ positive
coefficients of $U$ are therefore all $1$. Since $\varepsilon=1$,
we have $UV=Q$ and $\vnum_q(Q)\ge(r-1)+(p-r)=p-1$, contrary
to $\height(Q)\le p-2$.

\caseheading{2.2.3}{$l\ge2$}
We have $z_0^2\mid U$. For $d=p-2r+2$,
$$
 dg=\sigma(\overline{z_0^2})\in[-(r-1),r-1]g,\qquad
 0<d<p-r+1.
$$
Hence $d\le r-1$. The subsequence $U(z_0^2)^{-1}q^d$ of $A$
has image sum zero and length at most $2r-3<|A|$. By \eqref{eq:fixed-A}, its sum is $c=\sigma(U)$.
Thus $2z_0=dq$ and
$$
 2c=2z_0+2(r-1)q=pq.
$$
It follows that $(2P+ld)q=2\sigma(U)=pq$, whereas
$$
 p<p+(l-1)(p-2r)\le2P+l d\le(2r-l)(r-1)<p^2.
$$
This contradicts the order of $q$.
\end{proof}

\begin{proposition}\label{prop:double}
Let $S=a^{p-1}b^{p-1}xy$ satisfy \eqref{eq:common}, with $a\ne b$ and $x,y\notin\{a,b\}$. Then $S$ is exceptional if and only if there exists $i\in[1,p-2]$ such that
$$
 x=-ia-(p-1-i)b,\qquad y=-(p-1-i)a-ib.
$$
\end{proposition}
\begin{proof}
Since $\overline a\ne\overline b$, there are unique
$i,j\in[0,p-1]$ satisfying
$$
 i(\overline a-\overline b)=\overline b-\overline x,\qquad
 j(\overline a-\overline b)=\overline b-\overline y.
$$
Neither $\overline x$ nor $\overline y$ belongs to
$\{\overline a,\overline b\}$, so $i,j\in[1,p-2]$. Set $I=i+j$.
Let
$$
 T_1=a^ib^{p-1-i}x,\qquad T_2=a^jb^{p-1-j}y.
$$
The equations for $i$ and $j$ give $\sigma(T_1),\sigma(T_2)\in pG$.

\par\medskip
\noindent\textit{Necessity.}\enspace
Suppose that $S$ is exceptional. Since $|T_1|=|T_2|=p$,
we have $\sigma(T_1)=\sigma(T_2)=0$ by \eqref{eq:exceptional}. Hence
$$
 x=-ia-(p-1-i)b,\qquad y=-ja-(p-1-j)b.
$$
For $0\le u\le p-2$ we have
$$
 \sigma(a^ub^{p-2-u}xy)=(u-I)(a-b)-pb.
$$
If $I\ne p-1$, let $u$ be the least nonnegative residue of $I$
modulo $p$. Since $2\le I\le2p-4$, we have $u\in[0,p-2]$ and
$u-I\in\{0,-p\}$. The displayed sum is therefore either $-pb$ or
$-pa$, both in $pG\setminus\{0\}$. This contradicts exceptionality,
so $I=p-1$.

\par\medskip
\noindent\textit{Sufficiency.}\enspace
Suppose that the formulas in the proposition hold. By uniqueness
of the indices defined above, $j=p-1-i$ and $I=p-1$.
Let $T\mid S$ have length $p$, and let $u$ count its terms equal
to $a$. We separate the cases by the number of selected positions
among those of $x$ and $y$.

\caseheading{1}{Neither position is selected}
Then $1\le u\le p-1$ and $\sigma(\overline T)=u(\overline a-\overline b)\ne0$.

\caseheading{2}{Exactly one position is selected}
If $\sigma(\overline T)=0$, then $u=i$ for the position of $x$
and $u=j$ for the position of $y$. The formula for the selected
term gives $\sigma(T)=0$.

\caseheading{3}{Both positions are selected}
The equality $\sigma(\overline T)=0$ would require
$u\equiv p-1\pmod p$, which is impossible for $u\in[0,p-2]$.

\par Thus \eqref{eq:exceptional} holds, so $S$ is exceptional.
\end{proof}

\subsection{Other large multiplicities}

For $1\le t<m$, choose $I\subseteq[1,|R|]$ with $|I|=t$, and let
$$
B=\prod_{i\in I}(x_i-a),\qquad
R'=\prod_{i\in[1,|R|]\setminus I}x_i,\qquad
S'=a^{p-m+t}R'.
$$
Then $B\mid Q$ and
$$
\supp(S')\subseteq\supp(S),\qquad
\vnum_a(S')=p-m+t,\qquad
\vnum_x(S')\le \vnum_x(S)\le p-m \quad\text{for }x\ne a.
$$
Thus $S'$ satisfies \eqref{eq:common} whenever $S$ does, and $a$ is
the unique term of maximum multiplicity in $S'$.

\begin{lemma}\label{lem:replacement-high}
Suppose that $S$ is exceptional and that the sequence $S'$ obtained
by the above replacement is not exceptional. Then there exists
$W\mid QB^{-1}$ such that
\begin{equation}\label{eq:replace-high}
 m-t\le |W|<m,\qquad \sigma(\overline W)=0,\qquad \sigma(W)\ne c.
\end{equation}
\end{lemma}
\begin{proof}
Choose a subsequence $T\mid S'$ with $|T|=p$ and
$\sigma(T)\in pG\setminus\{0\}$. Subtracting $a$ from each term of
$T$ different from $a$ gives $W\mid QB^{-1}$, since these terms
come from unselected positions in $R$. The multiplicity of $a$
in $S'$ gives $|W|\ge m-t$. If $|W|\ge m$, the terms of $T$
different from $a$ can be combined in $S$ with $p-|W|\le p-m$
copies of $a$ to form $T$, a contradiction. Finally,
$\sigma(T)=pa+\sigma(W)\in pG\setminus\{0\}$ gives
$\sigma(\overline W)=0$ and $\sigma(W)\ne c$.
\end{proof}

\begin{proposition}\label{prop:23}
Under \eqref{eq:common}, if $m\in\{2,3\}$, then $S$ has a subsequence
$T$ with $|T|=p$ and $\sigma(T)\in pG\setminus\{0\}$.
\end{proposition}
\begin{proof}
Assume that $S$ is exceptional. We select $m-1$ positions in $R$
such that, for the corresponding subsequence $B\mid Q$,
$\overline{QB^{-1}}$ has no nonempty zero-sum subsequence of length
less than $m$.

\caseheading{1}{$m=2$}
Any position suffices because $0\notin\supp(\overline Q)$.

\caseheading{2}{$m=3$}
By \eqref{eq:common}, $|\supp(S)|\le p-1$, so
$2p\le(p-1)(p-3)$ and hence $p\ge7$. If $\overline Q$ had three
position-disjoint two-term zero-sum subsequences, denote their
corresponding subsequences of $Q$ by $B_1,B_2,B_3$. The pairwise
products have length $4$, and the product of all three has length
$6$. By \eqref{eq:fixed}, all four products have sum $c$. Hence
$$
 3c=\sigma(B_1B_2)+\sigma(B_1B_3)+\sigma(B_2B_3)
    =2\sigma(B_1B_2B_3)=2c,
$$
contrary to $c\ne0$. For each pair of opposite nonzero values in
$\overline Q$, select the positions in $R$ corresponding to all
occurrences of a value of smaller multiplicity. The number selected
in each pair equals the maximum number of disjoint two-term
zero-sum subsequences formed from that pair. Thus at most two
positions are selected in total. Deleting the corresponding
differences eliminates all two-term zero-sum subsequences of
$\overline Q$. Select further positions if necessary to make
exactly two.

\par In both cases, the number of selected positions is $t=m-1$.
Thus $a$ is the only term of $S'$ with multiplicity $p-1$. By Proposition~\ref{prop:base}, $S'$ is not
exceptional. By Lemma~\ref{lem:replacement-high}, there
is a nonempty $W\mid QB^{-1}$ of length less than $m$ with
$\sigma(\overline W)=0$, contrary to the choice of the positions.
\end{proof}

\begin{lemma}\label{lem:complement-high}
Suppose $S$ is exceptional, $m\ge4$, $p\ge3m-3$, and
$\sigma(D)\ne0$ for every $D\mid Q$ with $1\le|D|<m$.
Let $BW\mid Q$ with $1\le|B|,|W|<m$, $|BW|\ge m$, and
$\sigma(B),\sigma(W)\in pG$. Then
$$
 T\mid Q(BW)^{-1},\quad |T|\ge1,\quad \sigma(T)\in pG
 \quad\Longrightarrow\quad |T|\ge p-|B|+1.
$$
\end{lemma}

\begin{proof}
Suppose $|T|\le p-|B|$.

\caseheading{1}{$|T|<m$}
Since $m\le |BW|<|BWT|\le3m-3\le p$, \eqref{eq:fixed} gives
$$
 \sigma(T)=\sigma(BWT)-\sigma(BW)=0.
$$
This contradicts the nonzero-sum hypothesis for subsequences of length less than $m$.

\caseheading{2}{$m\le|T|\le p-|B|$}
Both $T$ and $TB$ have lengths in $[m,p]$ and sums in $pG$.
By \eqref{eq:fixed},
$$
 \sigma(B)=\sigma(TB)-\sigma(T)=0,
$$
contrary to the hypothesis on $B$.
\end{proof}

\begin{lemma}\label{lem:short-zero-sum-high}
Suppose $S$ is exceptional and satisfies \eqref{eq:common}, with
$m\ge4$ and $p\ge3m-3$. Assume that for every $B\mid Q$ with $t=|B|\in[1,m-1]$ there exists $W\mid QB^{-1}$
satisfying \eqref{eq:replace-high}. Then
$$
 \sigma(D)\ne0\quad\text{for every }D\mid Q\text{ with }1\le|D|<m.
$$
Moreover, if $B\mid Q$ is nonempty and $\sigma(B)\in pG$, then
$|B|>\lfloor(m+1)/2\rfloor$. In particular, if also $|B|<m$, then
$\overline B$ is minimal zero-sum.
\end{lemma}

\begin{proof}
To prove the first assertion, suppose that $1\le|B|<m$ and
$\sigma(B)\in pG$. Choose $W$ as in \eqref{eq:replace-high}.
Since $m\le|BW|\le2m-2\le p$, \eqref{eq:fixed} gives
$\sigma(B)+\sigma(W)=c$, so $\sigma(B)\ne0$.

Now suppose that $\sigma(B)\in pG$ and $t=|B|\in[1,\lfloor(m+1)/2\rfloor]$.
Choose $W$ as in \eqref{eq:replace-high} and let $V_0=Q(BW)^{-1}$.

\caseheading{1}{$\overline{V_0}$ is zero-sum free}
Take $T=V_0$. Then $|T|=p+m-t-|W|\ge p-t+1$.

\caseheading{2}{$\overline{V_0}$ has a nonempty zero-sum subsequence}
Choose $T\mid V_0$ with $\overline T$ minimal zero-sum.
Lemma~\ref{lem:complement-high} gives $|T|\ge p-t+1$.

\par In both cases, $p>2t$ gives $|T|>(p+2)/2$. By
Lemma~\ref{lem:long-normal-form}, the positive coefficients of
$\overline T$ relative to a suitable generator sum to at most $p$.
Lemma~\ref{lem:coefficient-multiplicity} gives
$$
 \height(\overline T)\ge2|T|-p\ge p-2t+2>p-m,
$$
contrary to $\height(\overline Q)\le p-m$. If $|B|<m$ and
$\overline B$ were not minimal zero-sum, its two nonempty zero-sum
factors would each have length at least $\lfloor(m+1)/2\rfloor+1$,
which is impossible.
\end{proof}

\begin{proposition}\label{thm:high}
Under \eqref{eq:common}, suppose that $m\ge4$ and
$p\ge\max\{3m-3,4m-11\}$. Then $S$ has a subsequence $T$ with
$|T|=p$ and $\sigma(T)\in pG\setminus\{0\}$.
\end{proposition}
\begin{proof}
We use strong induction on $m\ge4$ and assume that $S$ is exceptional.
For any choice of $t$ positions in $R$ with $1\le t<m$, let $B$
and $S'$ be as above. To show that $S'$ is not exceptional, we
distinguish three cases according to its parameter $m-t$.

\begin{itemize}[nosep]
\item If $m-t=1$, Proposition~\ref{prop:base} applies because $a$ is
the unique term of maximum multiplicity in $S'$.
\item If $m-t\in\{2,3\}$, Proposition~\ref{prop:23} applies.
\item If $m-t\ge4$, the induction hypothesis applies, since
$p\ge\max\{3(m-t)-3,4(m-t)-11\}$.
\end{itemize}

In all three cases, Lemma~\ref{lem:replacement-high}
gives $W\mid QB^{-1}$ satisfying \eqref{eq:replace-high}.

Thus the hypothesis of Lemma~\ref{lem:short-zero-sum-high} holds
for every choice of $a$ of maximum multiplicity. In particular,
every nonempty zero-sum subsequence of $\overline Q$ of length
less than $m$ is minimal zero-sum and has length at least three.
Also, no nonempty subsequence of $Q$ of length less than $m$
has sum zero in $G$.

Apply \eqref{eq:replace-high} with $t=1$ to obtain $U\mid Q$ with
$|U|=m-1$ and $\sigma(\overline U)=0$. Repeat the replacement step
with the positions in $R$ corresponding to $U$ to obtain
$W\mid QU^{-1}$ with $\sigma(\overline W)=0$. Let $Q=UWV$ and
$K=|UW|$. Then $m\le K\le2m-2\le p$ and $\sigma(UW)\in pG$,
so \eqref{eq:fixed} gives $\sigma(UW)=c$. In particular,
$\overline U$ and $\overline W$ are minimal zero-sum and have
length at least three.

By Lemma~\ref{lem:complement-high}, with $W$ in the role of $B$, every
nonempty zero-sum subsequence of $\overline V$ has length at least
$p-|W|+1=|V|$. Thus
$\overline V$ is zero-sum free or itself minimal zero-sum. Moreover,
$$
 |V|=p+m-K\ge p-m+2>(p+2)/2.
$$
By Lemma~\ref{lem:long-normal-form}, choose a generator $g$ of $C_p$
so that the terms of $\overline V$ have the form $sg$ with
$s\in[1,p-1]$ and these integer coefficients sum to $p-\delta$
for some $\delta\ge0$. Thus $\sigma(\overline V)=-\delta g$.
Since $2|V|>p$, Lemma~\ref{lem:coefficient-multiplicity} gives
an occurrence of $g$ in $\overline V$. By injectivity on $\supp(Q)$
there is a unique value $q$ in $\supp(V)$ with $\overline q=g$.
For $v=\vnum_q(V)$, the same lemma and $K\le2m-2$ give
\begin{equation}\label{eq:inventory-high}
 v\ge 2|V|-(p-\delta)=p+2m-2K+\delta\ge p-K+2.
\end{equation}
If at most one coefficient exceeds $1$, then
$v\ge |V|-1\ge p-m+1$, a contradiction. Write the coefficients as
$a_1,\ldots,a_{|V|}$. We have
$$
 \sum_{i=1}^{|V|}(a_i-1)=(p-\delta)-|V|=K-m-\delta.
$$
Since at least two summands are positive, each coefficient $s>1$
satisfies $2\le s\le K-m-\delta$.

Let $J=\{0\}\cup\Sigma(\overline{UW})$. By \eqref{eq:minsum-high} and Lemma~\ref{lem:cauchy-davenport}, $|J|\ge\min\{p,2K-3\}$. For each $x\in J\setminus\{0\}$ there is a subsequence $D\mid UW$
with $\sigma(\overline D)=x$ and $2\le|D|\le K-2$. For a one-term representation, adjoin the other complete zero-sum factor. For a representation of length $K-1$, remove the complete zero-sum factor it contains.

For each coefficient $s\ge2$ that occurs in $\overline V$, let
$w_s$ be the unique value in $\supp(V)$ with $\overline{w_s}=sg$.
Since $2\le s\le K-m-\delta<p$, we have $w_s\ne q$.
We claim that
\begin{equation}\label{eq:exchange-high}
 [m+\delta+s-3,p-K+\delta+2]g\cap J\ne\varnothing
\end{equation}
implies $w_s=sq$. Choose $z$ in the indicated interval with
$zg\in J$ and
$D\mid UW$ with $\sigma(\overline D)=zg$ and $2\le |D|\le K-2$.
The interval lies in $[1,p-1]$ since $\delta\le K-m-2$.
The image of the sum of $DV$ is $(z-\delta)g$, so we compare
the effect of removing either $w_s$ and $z-\delta-s$ copies
of $q$, or $z-\delta$ copies of $q$. Set
$$
 T_1=DV(w_sq^{z-\delta-s})^{-1},\qquad
 T_2=DV(q^{z-\delta})^{-1}.
$$
Both exponents are nonnegative and at most $p-K+2\le v$.
The sums of the two removed subsequences both have image
$(z-\delta)g$ under $\pi$, so
$$
 \sigma(\overline T_1)=\sigma(\overline T_2)
   =zg-\delta g-(z-\delta)g=0.
$$
For $d=|D|$ their lengths are
$$
\begin{aligned}
 |T_1|&=p+m-K-z+\delta+d+s-1,\\
 |T_2|&=p+m-K-z+\delta+d.
\end{aligned}
$$
The bounds on $z$ and $d$ give $m\le |T_2|\le |T_1|\le p$.
Hence \eqref{eq:fixed} gives $\sigma(T_1)=\sigma(T_2)=c$ and
$w_s=sq$.

We give a criterion for obtaining a contradiction from terms
satisfying $w_s=sq$. Let $\rho\in\{3,4\}$ and suppose that
$p\ge2K-2\rho$. Since $K\ge6>\rho$, the interval
$[K-\rho,p-K+\rho]$ is a nonempty subset of $[1,p-1]$.
The set $[K-\rho,p-K+\rho]g$ has complement in $C_p$ of size
$2K-2\rho-1<|J|$, since $|J|\ge\min\{p,2K-3\}$. Choose
$r$ in this interval with $rg\in J$ and choose $D\mid UW$ with
$\sigma(\overline D)=rg$ and $2\le|D|\le K-2$.

Suppose that $E\mid V$ consists of terms $w_s$ with $s>1$ and
$w_s=sq$. Let $k=|E|$ and let $n$ be the sum of the coefficients $s$
of its terms $w_s$. If
$$
 n\le K-\rho,\qquad \rho-2\le n-k\le K-\rho+2-m,
$$
then $DE$ and $(UW)D^{-1}E$ have sums with images $(r+n)g$
and $(n-r)g$, respectively. To cancel these images, consider
$$
 T_1=D E q^{p-r-n},\qquad
 T_2=(UW)D^{-1}E q^{r-n}.
$$
Both are subsequences of $Q$. Indeed, $n\le K-\rho$ ensures that
both exponents are nonnegative. Each is at most
$p-K+\rho-n\le p-K+2\le v$, since $n\ge n-k\ge\rho-2$.
Also, $E$ contains no term equal to $q$. The lengths satisfy
$$
\begin{aligned}
 |T_1|&=|D|+p-r-(n-k),\\
 |T_2|&=K-|D|+r-(n-k).
\end{aligned}
$$
Both lie in $[K-\rho+2-(n-k),p+\rho-2-(n-k)]\subseteq[m,p]$.
We have $\sigma(\overline{T_1})=\sigma(\overline{T_2})=0$, since the coefficients of the terms of
$E$ sum to $n$ and $\sigma(\overline D)=rg$.
Using $\sigma(E)=nq$ and \eqref{eq:fixed}, we get
$$
 2c=\sigma(T_1)+\sigma(T_2)=\sigma(UW)+pq=c+pq.
$$
Hence $pq=c=-pa$, contrary to $a+q$ being a unit.

We now distinguish two cases according to whether $p\ge2K-5$.

\caseheading{1}{$p\ge2K-5$}
Take $\rho=3$, so $p\ge2K-2\rho$. For every coefficient $s>1$,
the interval in \eqref{eq:exchange-high} is nonempty and its
complement has size $K+m+s-6\le2K-6<|J|$. Thus $w_s=sq$.
Take $E=w_s$ for any such coefficient. Then
$n=s\le K-m-\delta\le K-3$ and
$1\le n-k=s-1\le K-m-1$, so the criterion gives a contradiction.

\caseheading{2}{$p<2K-5$}
The bounds $p\ge4m-11$ and $K\le2m-2$, together with the parity of $p$,
give $p=2K-7=4m-11$ and $K=2m-2$. Since $p\ge3m-3$, we have
$m\ge8$. Take $\rho=4$, so $p\ge2K-2\rho$. Here $J=C_p$.
For each occurring coefficient $s\in[2,m-4]$, the interval
in \eqref{eq:exchange-high} has $m-s-3\ge1$ points, so $w_s=sq$.

\caseheading{2.1}{At least one coefficient lies in $[3,m-4]$}
Take $E=w_s$ for such a coefficient $s$. Then
$n=s\le K-4$ and $2\le n-k=s-1\le K-m-2$.
The criterion gives a contradiction.

\caseheading{2.2}{No coefficient lies in $[3,m-4]$ and at least two equal $2$}
Take $E=w_2^2$. Then $n=4\le K-4$ and
$n-k=2\le K-m-2$. The criterion gives a contradiction.

\caseheading{2.3}{No coefficient lies in $[3,m-4]$ and at most one equals $2$}
All other coefficients greater than $1$ are at least $m-3$.
At most one is that large, since $2(m-4)>m-2\ge K-m-\delta$.
Thus exactly two coefficients exceed $1$, one of them being $2$.
Now $v=|V|-2=p-m$. Since $q^v\mid Q$, the term $b=a+q$
occurs at least $v=p-m$ times in $S$ and thus has maximum
multiplicity. The two differences $a-b$ and $a+w_2-b$ satisfy
$$
 \overline{a-b}=-g,\qquad \overline{a+w_2-b}=g.
$$
Hence $a$ and $a+w_2$ are both different from $b$. The two
differences therefore occur in the sequence defined using $b$
in place of $a$. Their sum is $w_2-2q=0$ in $G$. This contradicts
Lemma~\ref{lem:short-zero-sum-high} applied with $b$ as the
term of maximum multiplicity.
\end{proof}

\section{Block interpolation and multiplicity bounds}\label{sec:foundation}

Throughout this section, we assume that $\pi$ is injective on
$\supp(S)$, as in \eqref{eq:common}.

Set $h=\height(S)=p-m$ and define
$$
 m_0=1+\min\!\left\{\left\lfloor\frac p3\right\rfloor+1,
                  \left\lfloor\frac{p+11}{4}\right\rfloor\right\}.
$$
For $p\ge23$, the condition $m\ge m_0$ is equivalent to
$m>(p+11)/4$. For $1\le j\le h$, define the $j$th layer by
$$
 L_j=\{x\in\supp(S):\vnum_x(S)\ge j\},\qquad s_j=|L_j|.
$$
Each value contributes one term to each layer it belongs to. In
particular, $\sum_{j=1}^h s_j=2p$. Set
$$
 f(s)=\left\lfloor\frac{s^2}{4}\right\rfloor,\qquad
 b_{S}=\sum_{j=1}^h f(s_j),\qquad e=2p-2-b_{S}.
$$
For a set $U$, write $\Sigma_k(U)$ for the sums of $k$ distinct
elements, with $\Sigma_0(U)=\{0\}$. A central selection from layer
$L_j$ takes $k_j\in\{\lfloor s_j/2\rfloor,\lceil s_j/2\rceil\}$
terms. The number of odd $s_j$ is even. Rounding equally often in
either direction on the odd layers gives
$$
 \sum_{j=1}^h k_j=\frac12\sum_{j=1}^h s_j=p.
$$
Specified odd layers may be paired and rounded in opposite directions
in advance.

\begin{lemma}\label{foundation:interpolation}
Let $L_1,\ldots,L_q\subseteq G$ be sets on each of which $\pi$ is
injective. For each $j\in[1,q]$, let $\ell_j$ be an integer with
$0\le \ell_j\le |L_j|$. Set
$$
 D=\sum_{j=1}^q \ell_j(|L_j|-\ell_j),\qquad Y=\sum_{j=1}^q\Sigma_{\ell_j}(L_j).
$$
Then $|\pi(Y)|\ge\min\{p,D+1\}$. If $Y\cap pG\subseteq\{0\}$, then $D\le2p-2$.
\end{lemma}

\begin{proof}
We give the same interpolation argument over $\mathbb Z/p\mathbb Z$
and $G=\mathbb Z/p^2\mathbb Z$. In this proof, $\cdot$ denotes
ring multiplication and $\prod$ denotes an indexed product in the
chosen coefficient ring or in the polynomial ring over it.

For each $j$, set $s_j=|L_j|$ and introduce variables
$X_{j,1},\ldots,X_{j,\ell_j}$, one for each position to be
filled from $L_j$.
A block with $\ell_j=0$ has no variables.
Integers occurring as polynomial coefficients are interpreted as
residue classes in the chosen ring. Since $\pi$ is injective on
$L_j$, we have $|\pi(L_j)|=s_j\le p$.

Fix an integer $d\in[0,D]$ and a monic polynomial $Q(T)$ of
degree $d$ over the chosen ring. Substituting
$\sum_{j=1}^q\sum_{i=1}^{\ell_j}X_{j,i}$ for $T$ in $Q(T)$ and
multiplying by the Vandermonde factors, we obtain
$$
 P=Q\!\left(\sum_{j=1}^q\sum_{i=1}^{\ell_j}X_{j,i}\right)\cdot
   \prod_{j=1}^q\prod_{1\le i<r\le\ell_j}(X_{j,r}-X_{j,i}).
$$
The Vandermonde factors make $P$ vanish whenever two variables
in the same block take the same value. If the choices are distinct
within each block, selecting them from $\pi(L_j)$ over
$\mathbb Z/p\mathbb Z$ gives a sum in $\pi(Y)$, and selecting them
from $L_j$ over $G$ gives a sum in $Y$.

We will detect the coefficient in $P$ of the monomial
$$
 X^{\boldsymbol t}=\prod_{j=1}^q\prod_{i=1}^{\ell_j}
                       X_{j,i}^{t_{j,i}}.
$$
For the variable $X_{j,i}$, the interpolation formula needs
$t_{j,i}+1$ distinct nodes from the $s_j$ available values in
block $j$, so we require
$0\le t_{j,i}\le s_j-1$. We choose the exponents within each
block to be strictly increasing. Their positive differences lie
in $[1,p-1]$, so the residue classes of those differences are units
in both rings. To make $X^{\boldsymbol t}$ have the top degree of
$P$, we also require
$$
 \sum_{j,i}t_{j,i}=d+\sum_{j=1}^q\binom{\ell_j}{2}.
$$
To see that these requirements can be met, set
$t_{j,i}=i-1+\lambda_{j,i}$. The inequalities within a nonempty
block become
$$
 0\le\lambda_{j,1}\le\cdots\le\lambda_{j,\ell_j}
 \le s_j-\ell_j.
$$
The total-degree requirement is now $\sum_{j,i}\lambda_{j,i}=d$.
For a block $j$ with $\ell_j>0$ and a desired sum
$d_j\in[0,\ell_j(s_j-\ell_j)]$ of the $\lambda_{j,i}$,
divide $d_j$ by $\ell_j$ with
quotient $a_j$ and remainder
$b_j$. Take $\ell_j-b_j$ entries equal to $a_j$ followed by
$b_j$ entries equal to $a_j+1$. These entries satisfy the
upper bound, since $a_j=s_j-\ell_j$ forces $b_j=0$.
A block with $\ell_j=0$ contributes zero. The intervals for
$d_j$ have sum $[0,D]$, so we may choose them with
$\sum_j d_j=d$. This gives the required exponents.

For each variable $X_{j,i}$, choose $E_{j,i}\subseteq L_j$ with
$|E_{j,i}|=t_{j,i}+1$. The node set for this variable is
$\pi(E_{j,i})$ over $\mathbb Z/p\mathbb Z$ and $E_{j,i}$ over $G$.
In either ring, let $E$ be the node set for $X_{j,i}$. For a
univariate polynomial $H(T)$ over that ring, define
$$
 \mathcal L_{j,i}(H)=\sum_{a\in E}
   \frac{H(a)}{\displaystyle\prod_{b\in E\setminus\{a\}}(a-b)}.
$$
Distinct nodes have distinct residues modulo $p$, so all node
differences in the denominators are units. The Vandermonde
determinant of the nodes is therefore a unit, and Lagrange
interpolation holds in either ring. Taking the coefficient of
$T^{t_{j,i}}$ in the interpolating expression for $T^u$ gives
$$
 \mathcal L_{j,i}(T^u)=
 \begin{cases}
 0,&0\le u<t_{j,i},\\
 1,&u=t_{j,i}.
 \end{cases}
$$
Every monomial in $P$ other than $X^{\boldsymbol t}$ has an
exponent below its target $t_{j,i}$ in at least one variable,
since $\deg P\le\sum_{j,i}t_{j,i}$. Applying the functionals
successively in their respective variables therefore extracts
the coefficient of $X^{\boldsymbol t}$. If $P$ vanishes on the
Cartesian product of the node sets in the chosen ring, this
coefficient is zero.

To compute the coefficient, only the leading term of $Q$ is
needed. In the determinant
$\det[X_{j,i}^{r-1}]_{i,r=1}^{\ell_j}$ for block $j$, using
column $r$ in row $i$ leaves the exponent $t_{j,i}-(r-1)$
to be supplied by the leading term of $Q$. The multinomial
expansion therefore gives
$$
 d!\cdot\prod_{j=1}^q
 \det\!\left[\frac{1}{(t_{j,i}-r+1)!}\right]_{i,r=1}^{\ell_j},
$$
where the entry is zero when $t_{j,i}<r-1$. All factorials
in the denominators are units, since $t_{j,i}<p$. Factoring
$1/t_{j,i}!$ from row $i$ of each determinant leaves the entry
$(t_{j,i})_{r-1}$, where $(t)_{r-1}=\prod_{k=0}^{r-2}(t-k)$
is the falling factorial and $(t)_0=1$. These are monic
polynomials of successive degrees $0,\ldots,\ell_j-1$,
so the determinants reduce to Vandermonde products. The
coefficient is thus
$$
 d!\cdot\prod_{j=1}^q
 \frac{\displaystyle\prod_{1\le i<r\le\ell_j}(t_{j,r}-t_{j,i})}
      {\displaystyle\prod_{i=1}^{\ell_j}t_{j,i}!}.
$$
An empty block contributes the factor $1$. Every factor other
than $d!$ is a unit. Therefore, if $d!$ is nonzero in the
chosen ring, then $P$ is nonzero at some tuple of nodes. At that
tuple the choices in each block are distinct, and $Q$ is nonzero
at their sum.

For the first assertion, work over $\mathbb Z/p\mathbb Z$ and set
$d=\min\{p-1,D\}$. If $|\pi(Y)|\le d$, the monic polynomial
$$
 Q(T)=T^{d-|\pi(Y)|}\cdot\prod_{a\in\pi(Y)}(T-a)
$$
has degree $d$ and vanishes on $\pi(Y)$, contrary to the
preceding conclusion because $d!$ is nonzero modulo $p$.
Hence $|\pi(Y)|\ge d+1=\min\{p,D+1\}$.

For the second assertion, work over $G$ and suppose that
$D\ge2p-1$ and $Y\cap pG\subseteq\{0\}$.
Take
$$
 d=2p-1,\qquad Q(T)=T\cdot\prod_{a=1}^{p-1}(T-a)^2.
$$
If $y\in Y$ has $\pi(y)=0$, then $y=0$. Otherwise let $a$ be the
integer representative of $\pi(y)$ in $[1,p-1]$. The residue class
of $a$ in $G$ satisfies $y-a\in pG$, so $(y-a)^2=0$ in $G$.
Thus $Q$ vanishes on $Y$. Since $p\mid d!$ but $p^2\nmid d!$,
the residue class of $d!$ in $G$ is nonzero. This again contradicts
the preceding conclusion.
\end{proof}

\begin{lemma}\label{foundation:budget}
Suppose $S$ is exceptional. Then $e\ge0$. Fix central selection
sizes $k_j$ of total $p$. For a set $J$ of layer indices, write
$$
 b=\sum_{j\in J}f(s_j),\qquad Y_J=\sum_{j\in J}\Sigma_{k_j}(L_j).
$$
We take $Y_{\varnothing}=\{0\}$. If $b\le p-e-1$, then $\pi$ is injective on $Y_J$ and $|\pi(Y_J)|\le b+e+1$.
\end{lemma}

\begin{proof}
Apply Lemma~\ref{foundation:interpolation} to the $h$ layers with
$\ell_j=k_j$ for $1\le j\le h$. For each layer, central selection
gives $k_j(s_j-k_j)=f(s_j)$. Thus
$$
 D=\sum_{j=1}^h k_j(s_j-k_j)=\sum_{j=1}^h f(s_j)=b_{S}.
$$
Also, $Y\subseteq\Sigma_p(S)$ because
$\sum_{j=1}^h k_j=p$. Thus $Y\cap pG\subseteq\{0\}$ by exceptionality,
and Lemma~\ref{foundation:interpolation} gives $b_{S}\le2p-2$.
Let $J^c$ consist of the remaining layer indices. Apply
Lemma~\ref{foundation:interpolation} to the layers indexed by
$J^c$ with $\ell_j=k_j$. The sumset is $Y_{J^c}$ and the parameter
is $D=b_{S}-b=2p-2-e-b$. Since $b\le p-e-1$, we have
$D\ge p-1$, so $\pi(Y_{J^c})=C_p$.
If $y,y'\in Y_J$ satisfy $\overline y=\overline{y'}$, choose
$z\in Y_{J^c}$ with $\overline z=-\overline y$. Then
$y+z,y'+z\in\Sigma_p(S)\cap pG$. By \eqref{eq:exceptional},
$y+z=y'+z=0$, so $y=y'$.

Apply Lemma~\ref{foundation:interpolation} again, taking one
element from $Y_J$ and $k_j$ elements from each layer indexed
by $J^c$. The resulting sumset is $Y_J+Y_{J^c}\subseteq\Sigma_p(S)$,
so exceptionality allows us to apply the second assertion of
Lemma~\ref{foundation:interpolation} with
$D=(|Y_J|-1)+(b_{S}-b)$. Thus
$$
 (|Y_J|-1)+(b_{S}-b)\le2p-2,\qquad
 |\pi(Y_J)|=|Y_J|\le b+e+1.
$$
\end{proof}

\begin{lemma}\label{foundation:short-filler}
If $S$ is exceptional and $1\le t\le\lfloor(m+1)/2\rfloor$, then $\pi$ is injective on $\Sigma_t(S)$.
\end{lemma}

\begin{proof}
Let $T_1,T_2\mid S$ with $|T_1|=|T_2|=t$ and
$\sigma(\overline{T_1})=\sigma(\overline{T_2})$. Outside their
combined positions choose $R$ of length $2p-2t$ and set
$H_0=p-2t+1$. Arrange equal terms consecutively and distribute the
positions cyclically into $H_0$ blocks. Since $2t\le m+1$,
we have $\height(R)\le p-m\le H_0$. Thus no block receives two
copies of the same value. Since $|R|-2H_0=2t-2\ge0$, each block has at least
two terms. Select half from each block, rounding up in half
of the odd-sized blocks and down in the other half, for a
total of $p-t$ terms. A block of size $u$ contributes $f(u)\ge u-1$,
so the total contribution is at least $|R|-H_0=p-1$.

Since the terms in each block are distinct and $\pi$ is injective
on $\supp(S)$, we may apply Lemma~\ref{foundation:interpolation}
to the sets of values in these blocks, with $\ell_j$ equal to the
number selected from block $j$. The parameter $D$ is at least
$p-1$, and the image under $\pi$ of the resulting sumset is
contained in $\Sigma_{p-t}(\overline R)$. Thus
$\Sigma_{p-t}(\overline R)=C_p$.
Choose $W\mid R$ with
$$
 |W|=p-t,\qquad
 \sigma(\overline W)=-\sigma(\overline{T_1})=-\sigma(\overline{T_2}).
$$
By \eqref{eq:exceptional}, $\sigma(T_1W)=\sigma(T_2W)=0$,
which gives $\sigma(T_1)=\sigma(T_2)$.
\end{proof}

Common completions also control the sums formed from the three most
frequent values of $S$.

\begin{lemma}\label{foundation:top-three}
Suppose $S$ is exceptional, $p\ge11$, and $m\ge m_0$. Let $v_1\ge v_2\ge v_3$ be the three largest multiplicities, padded with zeros if necessary. Then $v_2+v_3\le p-1$.
\end{lemma}

\begin{proof}
Suppose $v_2+v_3\ge p$. Then
$$
 0<m=p-v_1\le p-v_2\le v_3\le v_1=p-m.
$$
Thus $m\le(p-1)/2$. Choose distinct $a,b,c\in\supp(S)$ with
multiplicities $v_1,v_2,v_3$, respectively. The multiplicities of
any two of these values have sum at least $p$. Let $t$ be the maximum of
$\lfloor(m+2)/2\rfloor$ and the lower bound in the table. These
choices give $t<m$. The last row records the
upper bounds for the integer $M$ defined below.
$$
\begin{array}{c|rrrrrrr}
 p&11,13&17,19&23&29,31&37,41,43&47&\ge53\\ \hline
 t\text{ lower bound}&4&5&6&7&8&9&9\\
 M\text{ upper bound}&3&4&4&5&7&7&8
\end{array}
$$

For $(x,y)\in\mathcal T=\{(x,y)\in\mathbb Z_{\ge0}^2:x+y\le t\}$, define
$$
 T_{x,y}=a^x b^y c^{t-x-y},
 \qquad
 \sigma(T_{x,y})=tc+x(a-c)+y(b-c).
$$
Since each of $a,b,c$ occurs at least $v_3\ge m>t$ times, we have
$T_{x,y}\mid S$ for every $(x,y)\in\mathcal T$.

We first prove that $\pi$ is injective on
$$
 \{\sigma(T_{x,y}):(x,y)\in\mathcal T\}\subseteq G.
$$
Let $(x_i,y_i)\in\mathcal T$ for $i\in\{1,2\}$ and suppose that
$$
 \sigma(\overline{T_{x_1,y_1}})
 =\sigma(\overline{T_{x_2,y_2}}).
$$
If the two multiplicity vectors are equal, then the required equality
of sums is immediate. Hence assume that they are distinct. Then
$$
 x_1\ne x_2,\qquad y_1\ne y_2,
 \qquad x_1+y_1\ne x_2+y_2.
$$
Indeed, $x_1=x_2$ would give
$(y_1-y_2)(\overline b-\overline c)=0$, and hence $y_1=y_2$
because $|y_1-y_2|<p$ and $\overline b\ne\overline c$. The case
$y_1=y_2$ is analogous. If $x_1+y_1=x_2+y_2$, then
$(x_1-x_2)(\overline a-\overline b)=0$, which again forces the two
vectors to be equal.

Choose the copies in the two subsequences so that their positions are
nested within each of $a,b,c$. Their union then uses
$$
 \max\{x_1,x_2\},\quad
 \max\{y_1,y_2\},\quad
 \max\{t-x_1-y_1,t-x_2-y_2\}
$$
copies of $a,b,c$, respectively. Each of these three integers belongs
to $[1,t]$, so deleting the union leaves a copy of each of $a,b,c$.
Since $v_1\ge(p+1)/2$, the inequality $4v_1>|S|=2p$ shows that at
most three support values have multiplicity $v_1$. Hence every value
of multiplicity $v_1$ belongs to $\{a,b,c\}$.

Outside the union choose a subsequence $R$ satisfying
$$
 |R|=2p-2t,\qquad \supp(R)=\supp(S).
$$
Such a choice is possible because the remaining sequence has length
at least $2p-2t$ and
$2p-2t>p-1\ge|\supp(S)|$, where the last inequality follows from
the injectivity of $\pi$ on $\supp(S)\subseteq G^\times$. Every value
of multiplicity $v_1$ has
lost at least one copy, whereas every other value has multiplicity at
most $v_1-1$. Therefore
$$
 \height(R)\le v_1-1=p-m-1.
$$

Set $H_0=p-2t+1$. We construct a common completion according to the
maximum multiplicity of $R$.

\caseheading{1}{$\height(R)\le H_0$}
The blocks from the proof of Lemma~\ref{foundation:short-filler}
give $\Sigma_{p-t}(\overline R)=C_p$.

\caseheading{2}{$\height(R)>H_0$}
The inequalities
$$
 H_0<\height(R)\le p-m-1
$$
give $2t>m+2$, whereas the table gives $2t\le m+3$. Hence
$2t=m+3$, so $p-m-1=H_0+1$ and
$\height(R)=H_0+1$. Since $m\le(p-1)/2$, we have
$$
 4(H_0+1)-|R|=2p-3m-1\ge(p+1)/2>0.
$$

\caseheading{2.1}{$|\supp(R)|\ge4$}
First take a four-value block containing every value of multiplicity
$H_0+1$, of which there are at most three. The remainder has maximum
multiplicity at most $H_0$ and length $|R|-4=2H_0+2t-6\ge2H_0$.
Distribute it cyclically after grouping equal terms, obtaining $H_0$
blocks, each containing at least two terms with distinct values.
Select two terms from the four-value block and make central
selections from the remaining blocks, rounding the odd block sizes
in opposite directions in pairs. This selects $|R|/2=p-t$ terms.
The total contribution is at least $4+|R|-4-H_0=p-1$, so
Lemma~\ref{foundation:interpolation} gives $\Sigma_{p-t}(\overline R)=C_p$.

\caseheading{2.2}{$|\supp(R)|=3$}
Since $\{a,b,c\}\subseteq\supp(R)$, we have
$\supp(R)=\{a,b,c\}$. We distinguish whether their images form an
arithmetic progression.

\caseheading{2.2.1}{The three images form an arithmetic progression}
After permuting $a,b,c$ within this case, assume that
$\overline a+\overline c=2\overline b$.
Lemma~\ref{foundation:short-filler} with $t=2$ gives $a+c=2b$.
Set $\ell=\lceil(p-\vnum_b(S))/2\rceil$. Since any two of the three
multiplicities in $S$ sum to at least $p$, we have
$$
 \ell\le p-\vnum_b(S)
 \le\min\{\vnum_a(S),\vnum_c(S)\}.
$$
Also $p-2\ell\in\{\vnum_b(S)-1,\vnum_b(S)\}$. Consequently,
$$
 a^\ell c^\ell b^{p-2\ell}\mid S,
 \qquad
 \sigma(a^\ell c^\ell b^{p-2\ell})
 =\ell(a+c)+(p-2\ell)b=pb\in pG\setminus\{0\},
$$
contrary to \eqref{eq:exceptional}.

\caseheading{2.2.2}{The three images do not form an arithmetic progression}
Group equal terms of $R$ and distribute the positions cyclically
into $H_0+1$ blocks. Each block has distinct values since
$\height(R)=H_0+1$. The support has size three, and
$|R|-2(H_0+1)=2t-4\ge2$, so every block has size two or three
and at least two have size three. Select one term from one
three-value block and two from another. Their joint projected
sumset is a translate of
$$
 \{\overline a,\overline b,\overline c\}
 -\{\overline a,\overline b,\overline c\},
$$
which has seven elements since the three images do not form an
arithmetic progression. The two blocks supply three terms.
Round the remaining odd block sizes in opposite directions in
pairs to select $(|R|-6)/2$ more terms. Each remaining block has a projected
central sumset of cardinality equal to its size. Applying
Lemma~\ref{lem:cauchy-davenport} to the joint sumset and these
remaining sumsets gives
$$
 |\Sigma_{p-t}(\overline R)|
 \ge\min\{p,1+6+(|R|-6)-(H_0-1)\}=p.
$$

Cases~1, 2.1, and 2.2.2 give $\Sigma_{p-t}(\overline R)=C_p$, so
the two subsequences $T_{x_1,y_1}$ and $T_{x_2,y_2}$ have a common
completion $W\mid R$ satisfying
$$
 |W|=p-t,\qquad
 \sigma(\overline W)=-\sigma(\overline{T_{x_1,y_1}})
 =-\sigma(\overline{T_{x_2,y_2}}).
$$
By \eqref{eq:exceptional},
$$
 \sigma(T_{x_1,y_1}W)=\sigma(T_{x_2,y_2}W)=0,
$$
and hence $\sigma(T_{x_1,y_1})=\sigma(T_{x_2,y_2})$. Case~2.2.1
has already given a contradiction. This proves the required
injectivity of $\pi$.

We next count the sums $\sigma(T_{x,y})$. Since $\pi$ is injective on
their set, there are at most $p$ distinct sums in $G$. The table
gives $|\mathcal T|=\binom{t+2}{2}>p$ for $p<53$. For $p\ge53$,
$t\ge\lfloor(p+23)/8\rfloor$, so $p\le8t-17$. For $t\ge11$ this
gives $p<9t-27$, and the same inequality follows from $p\le53$
when $t=9$ and $p\le61$ when $t=10$. Since $t\ge9$, we have
$\binom{t+2}{2}-(9t-27)=(t-7)(t-8)/2>0$. Thus two distinct points
of $\mathcal T$ give the same sum in $G$.

For distinct $(x_1,y_1),(x_2,y_2)\in\mathcal T$ giving the same sum,
call
$$
 (r,s)=(x_1-x_2,y_1-y_2)
$$
a collision vector. Every collision vector satisfies
$r(a-c)+s(b-c)=0$. Hence, for two collision vectors
$(r_1,s_1)$ and $(r_2,s_2)$ we have $r_i(a-c)+s_i(b-c)=0$ for
$i=1,2$.
Eliminating $b-c$ gives
$$
 (r_1s_2-r_2s_1)(a-c)=0,\qquad
 |r_1s_2-r_2s_1|\le2t^2<p^2.
$$
Since $a-c$ has additive order $p^2$, the determinant is divisible
by $p^2$ and hence is zero. All collision vectors are therefore collinear.
The greatest common divisor of the coordinates of a nonzero
collision vector is at most $t<p$, so multiplication by it is an
automorphism of $G$. Cancelling it gives a primitive relation
$(r,s)$. Set
$M=\max\{|r|,|s|,|r+s|\}\le t$.
Each class of equal sums in $\mathcal T$ is an integer chain
with step $(r,s)$, and
$\mathcal T\cap(\mathcal T+(r,s))$ is a translate of the integer
triangle with side parameter $t-M$. Counting points with no
predecessor gives the number of distinct sums as
$$
 \binom{t+2}{2}-\binom{t-M+2}{2}
 =\frac{M(2t+3-M)}2\le p.
$$
The expression is strictly increasing for $0\le M\le t$.
Substituting the first excluded value gives the bounds in the
table. For $p\ge53$, the value at $M=9$ is $9t-27>p$.

The differences of distinct values have additive order $p^2$.
Thus none of $r,s,r+s$ can be zero, since a zero coordinate or
$r+s=0$ would give a nonzero multiple of a single support difference
with coefficient of absolute value at most $t<p$. After permuting
$a,b,c$, we may assume that $ua+vc=Mb$, where $u+v=M$,
$1\le u\le v$, and $\gcd(u,v)=1$. Let $H,L$ be the multiplicities
of $a,c$, respectively, and put $N=p-\vnum_b(S)$ and
$j=\lceil N/M\rceil$. Then $H,L\ge N\ge m$. The table gives
$m\ge M-1$, so $ju\le(N+M-1)/2\le N$ and
$N\le jM\le N+M-1\le p$. If $jv\le L$, then
$a^{ju}c^{jv}b^{p-jM}\mid S$ and
$$
 \sigma(a^{ju}c^{jv}b^{p-jM})
 =j(ua+vc)+(p-jM)b=pb\in pG\setminus\{0\},
$$
a contradiction.

It remains to consider $jv>L$. The inequalities
$$
 (j-1)M+1\le N\le L\le jv-1
$$
give $ju\le M-2$. We must have $u=1$. Indeed, $u\ge2$ implies
$M\ge5$, hence $p\ge29$ and $N\ge11$. Since $M\le8$ and $u,v$ are
coprime, $v\le5$. Also $ju\le6$ forces $j\le2$, since $j\ge3$
would require $u=2,M=8,v=6$. Thus $jv\le10<N$, a contradiction.

Now $j\le v-1$ gives $N\le v(v-1)-1$. For $v\le2$ this
contradicts $N\ge5$. For $v=3$, $M=4$ gives $p\ge17$ and $N\ge7$,
again impossible. For $v=4$, $M=5$ gives $p\ge29$ and $N\le11$
allows only $p=29,31$. Then $j=3$ and $N=L=11$. In the remaining
cases $p\ge37$ and $5\le v\le7$.

Since $\overline c\ne\overline b$, the element $d=c-b$ is a unit.
The relation $a+vc=(v+1)b$ gives $a=b-vd$. Reserve one term equal to
$c$, choose $x$ terms equal to $a$ and $y$ terms equal to $c$,
and complete to $p-1$ terms with $b$. For
$x\in[x_0,H]$, where $x_0=\max\{0,v-(L-N)\}$, the admissible
values of $y$ form $[y_-(x),y_+(x)]$, with
$$
 y_-(x)=\max\{0,N-1-x\},\qquad
 y_+(x)=\min\{L-1,p-1-x\}.
$$
For $x\le N-1$ its width is $\min\{L-N+x,p-N\}\ge v$.
For $x\ge N-1$ the width is at least $m-1\ge v$.
Restoring the reserved term gives
$$
 \sigma(a^x c^{y+1} b^{p-1-x-y})=pb+(1-vx+y)d.
$$
As $x$ increases by one, the left endpoint of the coefficient
interval $[1-vx+y_-(x),1-vx+y_+(x)]$ decreases by $v$ or $v+1$.
The intervals have width at least $v$, so consecutive intervals
overlap or are adjacent, and their union is an integer interval. Since $x_0\le v<m\le p-L$ and $H\ge N$,
the selections $(x_0,L-1)$ and $(H,0)$ are admissible.
The width is therefore at least
$$
 vH+L-1-vx_0\ge p-1+(v-1)H-vx_0.
$$
For $v=4$, $H\ge p-11$ and $x_0=4$, so $(v-1)H-vx_0\ge3p-49\ge p$. For $5\le v\le7$, if $x_0=0$ then $(v-1)H\ge4m>p$. Otherwise $L=N+v-x_0\le H+v-x_0$, and $H+L\ge p$ gives
$$
 (v-1)H-vx_0\ge \frac{v-1}{2}p-v^2\ge p
 \qquad (p\ge37).
$$
The interval has width at least $2p-1$, so it contains $kp$ and
$(k+1)p$ for some integer $k$. The corresponding sums belong to
$\Sigma_p(S)\cap pG$ and differ by $pd\ne0$, contrary to
\eqref{eq:exceptional}.
\end{proof}

For the layers of $S$, list the multiplicities $v_i$ in decreasing
order, padding with zeros. In the expression
$b_{S}=\sum_i\lfloor i/2\rfloor v_i$, removing the four most frequent
values decreases the weight of each $v_i$ with $i\ge4$ by two
and removes the contributions of $v_2$ and $v_3$. Since
$v_1=p-m$ and $\sum_i v_i=2p$, it follows that
\begin{equation}\label{eq:tail-weight-identity}
 b_{S}-\sum_j f(\max\{0,s_j-4\})
 =v_2+v_3+2\sum_{i\ge4}v_i=2(p+m)-(v_2+v_3).
\end{equation}

\begin{lemma}\label{foundation:tail}
Under the hypotheses of Lemma~\ref{foundation:top-three}, set $s=s_1$. Then $s\ge4$ and
$$
 f(s-4)\le p-2m-e-3,\qquad f(s)\le p-2m-e+2(s-4)+1.
$$
Assign weights $4,8,7,6$ to layers of sizes $2,3,4,5$, respectively, other than the first layer. Let $\Phi$ be their total weight. Then $\Phi\ge-p+9m+2+5e$. If $p\ge23$, then $\Phi>(5p+107)/4+5e$.
\end{lemma}

\begin{proof}
If $s\le3$, then $v_2+v_3=2p-v_1\ge p$, contradicting
Lemma~\ref{foundation:top-three}.
By Lemma~\ref{foundation:top-three}, $v_2+v_3\le p-1$. Using this in
\eqref{eq:tail-weight-identity} and keeping only the first layer
in the sum gives the first inequality. The second follows from $f(s)-f(s-4)=2s-4$.

Assign a layer of size $u$ the weight $w(u)=9(u-1)-5f(u)$.
This is zero for $u=1$, takes the four stated values for
$u=2,3,4,5$, and is nonpositive for $u\ge6$. Since
$\sum_j(s_j-1)=p+m$, the total weight is
$$
 \sum_j w(s_j)=9(p+m)-5(2p-2-e)=-p+9m+10+5e.
$$
Removing the first layer costs at most eight. Discarding the remaining
nonpositive weights gives $\Phi\ge-p+9m+2+5e$.
For $p\ge23$, substituting $m>(p+11)/4$ gives
$\Phi>(5p+107)/4+5e$.
\end{proof}

\section{Sumset growth and intervals from small layers}\label{sec:local-groups}

Retain the notation of Section~\ref{sec:foundation} and assume
\eqref{eq:common}. Suppose $p\ge23$, $m>(p+11)/4$, and $S$ is
exceptional with $e\ge1$.
Lemma~\ref{foundation:tail} gives $e\le(p-19)/2$. Consider the
nonfirst layers of sizes two through five. Layers of equal size have
the same support. For each odd size, pair the available layers and
round in opposite directions in each pair. We call such a pair a
balanced pair. Balance the remaining odd layers to select $p$ terms
in total.

For a finite support $U$ and $k\in[0,|U|]$, complementation gives
\begin{equation}\label{eq:sum-complement}
 \Sigma_{|U|-k}(U)=\left\{\sum_{x\in U}x\right\}-\Sigma_k(U).
\end{equation}
Thus the joint central sumset of two paired layers of size three
is a translate of $U-U$. For two paired layers of size five it
is a translate of $X-X$, where $X=\Sigma_2(U)$.

An integer model of a support $U\subseteq G$ has the form
$U=\{a+xd:x\in U_{\mathbb Z}\}$, where $d$ is a unit and the
coordinates in $U_{\mathbb Z}\subseteq\mathbb Z$ are distinct
modulo $p^2$. For a nonempty integer set $Y$, let
$\operatorname{diam}(Y)=\max Y-\min Y$.
A group of layers provides an initial interval if its central
sumset in $G$ contains an interval in integer coordinates with a
unit common difference. The width of that interval is its
diameter. The truncation loss is the diameter of the full central
sumset in those coordinates minus the width of the interval.
The constructions below also give integer coordinates for the
common support of the group.

\begin{lemma}\label{lem:local-integer-model}
Let $U\subseteq G$ have size four or five and set $A=\Sigma_2(U)$.
If $|U|=4$, assume $p\ge37$ and that $\pi$ is injective on $4A$.
If $|U|=5$, assume $p\ge29$ and that $\pi$ is injective on $A-A$.
Consider the expressions for elements of $4A$ in the first case
and of $A-A$ in the second. There is a rational affine model of
$U$ in which two such expressions are equal exactly when their
values in $G$ have the same image under $\pi$.
If the model has affine dimension one, then $U$ admits primitive
integer coordinates with a unit common difference and with the
same equalities.
\end{lemma}
\begin{proof}
In both cases a translate of $A$ lies in the set on which $\pi$
is injective. For two distinct elements of $U$, adjoining the
same third element gives distinct elements of $A$. Thus $\pi$
is injective on $U$ and all differences of distinct elements
of $U$ are units. Set $k=|U|$ and let $U=\{u_1,\ldots,u_k\}$.
Let $W=\{0,e_1,\ldots,e_{k-1}\}$, where the $e_i$ are the
standard basis vectors of $\mathbb Z^{k-1}$. Define the additive
homomorphism $\varphi\colon\mathbb Z^{k-1}\to G$ by
$\varphi(e_i)=u_{i+1}-u_1$. For $B=\Sigma_2(W)$ we have
$\varphi(B)=A-2u_1$. Take $Z=4B$ when $k=4$ and $Z=B-B$ when $k=5$.
Thus $\varphi(Z)=4A-8u_1$ or $A-A$, respectively. Let
$$
 L=\operatorname{span}_{\mathbb Q}
   \{v\in Z-Z:\pi(\varphi(v))=0\}.
$$
For each generating relation $v=z-z'$ with $z,z'\in Z$,
injectivity on $\varphi(Z)$ gives $\varphi(z)=\varphi(z')$. Thus
$\varphi(v)=0$. For the norm estimates, append
$-\sum_{i=1}^{k-1}v_i$ to each relation vector $v$ so that its
$k$ coordinates sum to zero. For $k=4$, each absolute coordinate is at most $4$
and the sums of the positive and absolute negative coordinates
are at most $8$, so the Euclidean norm is at most $8$. For $k=5$
the respective bounds are $2$, $4$, and $4$.

Let $M$ be an integer matrix with $k-1$ such relations as rows.
Reducing its entries modulo $p^2$, we see that $M$ annihilates the
column vector of support differences in $G$. The adjugate identity
then gives $(\det M)\cdot(u_i-u_1)=0$ in $G$ for $2\le i\le k$.
Each $u_i-u_1$ is a unit, so $p^2\mid\det M$. Hadamard's inequality gives
$$
 |\det M|\le
 \begin{cases}
 8^3,&k=4,\\
 4^4,&k=5,
 \end{cases}
 \qquad |\det M|<p^2.
$$
Thus $\dim L\le k-2$. If $L\ne\{0\}$, choose a basis from the
relation vectors.
For $k=4$, the minors of order at most two have absolute value
at most $32<p$. For $k=5$, a minor of order three has absolute
value at most $32$, since the maximum over $[-2,2]^9$ is attained
at a vertex and a $3\times3$ sign matrix has determinant divisible
by $4$ and of absolute value less than $8$. If such a minor is
odd, each row has an odd entry and norm at most $3$, so its
absolute value is at most $27<p$. Minors of orders one and two
have absolute value at most $8$. For $k=5$, a nonzero minor of
absolute value at most $32<2p$ divisible by $p$ would have to be
$\pm p$, contrary to the bound for odd minors. Hence a nonzero
minor of full rank is a unit modulo $p^2$.

If $L=\{0\}$, then every $v\in(Z-Z)\cap L$ evaluates to zero.
Otherwise, multiplying such a vector by the determinant of the
chosen minor clears the denominators in its expansion in the
basis. Each basis relation vanishes under $\varphi$, and the
determinant is a unit modulo $p^2$, so $\varphi(v)=0$.
Conversely, any $v\in Z-Z$ with $\pi(\varphi(v))=0$ belongs to
$L$ by definition. Thus for all $z,z'\in Z$,
$$
 z-z'\in L
 \quad\Longleftrightarrow\quad
 \pi(\varphi(z))=\pi(\varphi(z')).
$$
Thus the quotient $\mathbb Q^{k-1}/L$ gives the required model.
Translations into $Z$ also preserve the equalities in $B$, and
in $2B$ when $k=4$. If the quotient has dimension one, the
integer kernel of the matrix whose rows are the chosen basis
relations is generated by a primitive vector $n$. The free coordinate of $n$ is, up to sign, the
invertible minor divided by the greatest common divisor of the
maximal minors, so it is a unit modulo $p^2$. Elimination over
$G=\mathbb Z/p^2\mathbb Z$ then gives
$$
 (u_2-u_1,\ldots,u_k-u_1)=d n,
$$
where $d\in G$ multiplies each integer coordinate of $n$. The
free coordinate of the support-difference vector is a unit, so
$d$ is a unit. This gives the required integer coordinates.
\end{proof}

For any finite set $U$ and $1\le k<|U|$, we have
\begin{equation}\label{eq:support-difference}
 U-U\subseteq\Sigma_k(U)-\Sigma_k(U).
\end{equation}
Indeed, to express $u-u'$ for $u\ne u'$, adjoin the same $k-1$
elements of $U\setminus\{u,u'\}$ to $u$ and $u'$.
In primitive integer coordinates, take $Y=\Sigma_2(U_{\mathbb Z})$
or $Y=2\Sigma_2(U_{\mathbb Z})$. In both cases the differences in
$Y-Y$ have greatest common divisor $1$, since they include the
support differences by \eqref{eq:support-difference}. For the
twofold sumset, this inclusion follows by adding a common
summand to both terms of each difference. Translate
$Y$ and $-Y$ separately to have minimum zero. They have equal
diameters and cardinalities, so the parameter in the following
bound of Lev and Smeliansky~\cite[Lemma~1]{LevSme} is $\delta=1$.
\begin{equation}\label{eq:local-diameter}
 |Y-Y|\ge |Y|+\min\{\operatorname{diam}(Y),2|Y|-3\}.
\end{equation}
If the rational model has dimension at least two, then either
of the corresponding sumsets $Y$ in that model satisfies
\begin{equation}\label{eq:planar-difference}
 |Y-Y|\ge3|Y|-3.
\end{equation}
By \eqref{eq:support-difference}, the support differences belong
to $Y-Y$. A generic linear projection to the plane preserves the
cardinality of $Y$ and gives a two-dimensional set, so
Ruzsa's bound~\cite[Corollary~1.1]{RuzsaDim} gives
\eqref{eq:planar-difference}.

\begin{lemma}\label{lem:small-layer-local}
Either one of the constructions below gives an initial interval, or none applies. If an interval is available, choose one with the largest support size $q\in[2,5]$. Its width is at least both the sum of the contributions $f(s_j)$ of its layers and the integer diameter of its support. It uses at most four layers, whose total contribution is at most $12$. The truncation loss is at most $12$, and at most $8$ if $p\le31$. Every layer other than the first of size in $[2,q]$ can be written in integer coordinates with the same unit common difference. The gaps between successive central sums in such a layer are at most the width of the initial interval.

Call all nonfirst layers of sizes in $[2,q]$ aligned, including those used for the initial interval. If no construction applies, there are no aligned layers. The remaining nonfirst layers of sizes two through five can be partitioned into full groups and leftover layers as follows. For each full group with layer index set $J$, we have
$|\pi(Y_J)|\ge b+g+1$, where $b=\sum_{j\in J}f(s_j)$ and $g$ is
given in the table below.
$$
\begin{array}{c|c|r|r|c}
\text{Layer size}&\text{Layers per group}&b&g&\text{Maximum leftovers}\\ \hline
3&4&8&10&3\\
4&4&16&16&3\\
5&2&12&12&1
\end{array}
$$
No unaligned layer has size two. Each full group satisfies $b\le g$ and has total weight at most $16g/5$. The leftover layers have total weight at most $51$.
\end{lemma}

\begin{proof}
We first record how to lift progressions and their differences.
If $\pi$ is injective on $A+A$ for $A\subseteq G$, then
$\pi(a-a')=\pi(b-b')$ implies $\pi(a+b')=\pi(b+a')$, so
$a-a'=b-b'$. Thus injectivity on $A+A$ implies injectivity on
$A-A$. If $\pi$ is injective on $A-A$, it is injective on $A$
as well. Hence an arithmetic progression $\pi(A)$ lifts to an
arithmetic progression $A$ by equality of successive differences.

If $A=\Sigma_2(U)$ has integer coordinates of diameter $L\le16$
with unit common difference $d$, fix $u_0\in U$ and choose
$x_u\in[-L,L]$ with $u-u_0=x_ud$ for each $u\in U$, using
\eqref{eq:support-difference}. Since $u-v\in A-A$, the difference
$x_u-x_v$ is congruent modulo $p^2$ to an integer in $[-L,L]$.
As $3L<p^2$, these integers are equal. Thus the coordinates $x_u$
give $U$ diameter at most $L$.

The groups below use at most four layers of size four, four of
size three, or two of size five. Their total contribution is at
most $16$, and $p-e-1\ge(p+17)/2\ge20$. Hence
Lemma~\ref{foundation:budget} gives injectivity of $\pi$ on their
joint central sumsets.

We now consider each layer size separately.

\caseheading{1}{Layers of size two}
A two-element layer gives an interval of width $1$.

\caseheading{2}{Layers of size three}
Let $U$ be the three-element support and let
$\pi(U)=\{a_1,a_2,a_3\}$. To distinguish the supports that give
an interval from those that give growth, consider
$2\pi(U)-2\pi(U)$. The difference
$(c_1,c_2,c_3)$ of two distinct coefficient vectors with the same
sum in $C_p$ satisfies
$$
 \sum_{i=1}^3c_i=0,\qquad
 \sum_{c_i>0}c_i\le4,\qquad \sum_{i=1}^3 c_i a_i=0.
$$
If one coordinate were zero, the other two would give a nonzero
integer multiple of some $a_i-a_j$ equal to zero, with coefficient
of absolute value at most $4<p$. Thus all three coordinates are
nonzero. Up to permutation and sign the possibilities are
$$
 (1,1,-2),\quad(1,2,-3),\quad(1,3,-4),\quad(2,2,-4).
$$
Thus a collision requires $\pi(U)$ to be a translate of $\{0,g,rg\}$
in $C_p$ for some $g\ne0$ and $r\in\{2,3,4\}$.

\caseheading{2.1}{$\pi(U)$ is a translate of $\{0,g,rg\}$ for some $g\ne0$ and $r\in\{2,3,4\}$}
Let $u_0,u_1,u_r$ denote the support values
whose images have coordinates $0,1,r$ in this normalization. The relation
$$
 r\overline u_1=(r-1)\overline u_0+\overline u_r
$$
lifts to $G$ by Lemma~\ref{foundation:short-filler} whenever at
least $r$ copies of each value are available, counting the first
layer. Indeed, $r\le4\le\lfloor(m+1)/2\rfloor$. The lifted
relation gives $u_r=u_0+r(u_1-u_0)$, so the support has integer
coordinates $\{0,1,r\}$ with common difference $u_1-u_0$. For $r=2$
one layer gives $[0,2]$ up to translation. For $r=3$ one balanced
pair gives a translate of $[-3,3]$ by \eqref{eq:sum-complement}.
The table records the layer counts and widths.
$$
\begin{array}{c|ccc}
 r&2&3&4\\ \hline
 \text{Layers used}&1&2&4\\
 \text{Width}&2&6&16\\
 \text{Truncation loss}&0&0&0
\end{array}
$$
For $r=4$, the width follows from
$$
 (\{0,1,4\}-\{0,1,4\})+(\{0,1,4\}-\{0,1,4\})=[-8,8],
$$
since each difference set is $\{0,\pm1,\pm3,\pm4\}$.
If there are too few layers for the construction, at most three
remain.

\caseheading{2.2}{All other three-element supports}
A full group in this case needs four nonfirst layers of size three.
When four are available, use two balanced pairs.
The image under $\pi$ of their joint central sumset is a translate
of $2\pi(U)-2\pi(U)$. Its coefficient vectors have coordinate
sum zero and positive coordinates summing
to at most $2$. The counts for positive sums $0,1,2$ are $1,6,12$.
The classification above excludes collisions, so
$|2\pi(U)-2\pi(U)|=1+6+12=19$.

\caseheading{3}{Layers of size four}
If there are fewer than four nonfirst layers of size four, leave
them as leftovers. Otherwise let $U$ be their common support
and set $A=\Sigma_2(U)$ and $X=\pi(A)$. Two distinct pairs with equal sums
under $\pi$ are disjoint, since the images of the support values
are distinct. Two distinct complementary pairing equalities
would force two support images to coincide because $p$ is odd.
Thus $|X|=5$ or $6$.
Complementary pairs show that $A$ and $X$ are centrally symmetric.
Since a translate of $A$ is contained in $4A$, $\pi$ is injective
on $A$. Centering $X$ shows that $4X$ is symmetric and contains
its center, so $|4X|$ is odd. If $|4X|\ge33$, the four layers give
the desired growth. We may therefore assume $|4X|\le31$.

\caseheading{3.1}{$|X|=5$}
The set $A$ also has five points and contains its center. Translate
it to $A=\{0,\pm u,\pm v\}$. The $41$ integer pairs with
$|i|+|j|\le4$ give at most $31$ distinct values $i\overline u+j\overline v$
in $C_p$. Subtracting two pairs that give the same value yields a
nonzero relation $i\overline u+j\overline v=0$ with $|i|+|j|\le8$.
Both coefficients are nonzero because $\overline u,\overline v\ne0$
and $8<p$. Injectivity on $4A$ lifts the equality of the two
values to $G$, giving $iu+jv=0$. Cancelling the common divisor of
$i,j$ and interchanging or negating generators as needed gives
$u=ad,v=bd$ with a unit $d$, $1\le a<b$, $\gcd(a,b)=1$, and
$a+b\le8$. The integer fourfold sumset is
$\{ai+bj:|i|+|j|\le4\}$. Its diameter is at most $56<p^2$, so
multiplication by $d$ maps it injectively into $G$. By injectivity
of $\pi$ on $4A$, its cardinality is $|4X|$. The cases with
$|4X|\le31$ are listed below.
$$
\begin{array}{c|r|c|r|r}
(a,b)&|4X|&\text{Initial interval}&\text{Layers}&\text{Loss}\\ \hline
(1,2)&17&[-2,2]&1&0\\
(1,3)&23&[-4,4]&2&4\\
(1,4)&29&[-5,5]&2&6\\
(2,3)&25&[-6,6]&2&0
\end{array}
$$
The other ratios are $(1,b)$ for $5\le b\le7$, and $(2,5),(3,4),(3,5)$. The same formula gives $33,37,39,35,33,37$ points, respectively. Each interval in the table has width at least the diameter of $A$, hence at least that of $U$.

\caseheading{3.2}{$|X|=6$}
If $p\le31$, then
Lemma~\ref{foundation:budget} gives $|4X|\le17+e$, whereas
$$
 |2X|\ge13\quad\Longrightarrow\quad
 |4X|\ge\min\{p,25\}>17+e.
$$
Thus $|2X|=11$, since Lemma~\ref{lem:cauchy-davenport} gives $|2X|\ge11$
and symmetry makes $|2X|$ odd. For $p\ge37$,
Lemma~\ref{lem:cauchy-davenport} and $|4X|\le31$ give
$11\le |2X|\le16$, so $|2X|\in\{11,13,15\}$.

\caseheading{3.2.1}{$|2X|=11$}
Vosper's theorem in the form \eqref{eq:vosper} makes $X$ a
six-term arithmetic progression, which lifts to a one-layer interval.

\caseheading{3.2.2}{$|2X|\in\{13,15\}$}
Here $p\ge37$. Lemma~\ref{lem:local-integer-model} gives a rational model
preserving the cardinalities of $X$, $2X$, and $4X$.
By \eqref{eq:sum-complement}, the models of $X$ and $2X$ are
centrally symmetric, so their twofold sumsets are translates of
their difference sets. If the model has dimension at least two,
\eqref{eq:planar-difference} gives $|2X|\ge3|X|-3=15$, excluding
$|2X|=13$. Applying \eqref{eq:planar-difference} to the model of
$2X$ then gives $|4X|\ge3|2X|-3=42$, again a contradiction.
Thus the model has dimension one. Let $A_{\mathbb Z}=\Sigma_2(U_{\mathbb Z})$
in the primitive integer coordinates of Lemma~\ref{lem:local-integer-model}.
Then $|A_{\mathbb Z}|=|X|$, $|2A_{\mathbb Z}|=|2X|$, and
$|4A_{\mathbb Z}|=|4X|$.

If $|2X|=13$, then \eqref{eq:local-diameter} bounds the diameter of
$A_{\mathbb Z}$ by $7$. Apart from an arithmetic progression, a
six-point symmetric set of diameter at most $7$ is, up to
translation, $[0,6]\setminus\{3\}$ or is obtained from $[0,7]$ by
deleting a symmetric pair. Deleting $\{2,5\}$ or $\{3,4\}$ gives a
twofold sumset of size $15$, leaving only $[0,7]\setminus\{1,6\}$.
The two remaining shapes give intervals $[0,12]$ and $[2,12]$ using
two layers, with losses $0$ and $4$.

If $|2X|=15$, the diameter of $2A_{\mathbb Z}$ is at least $14$
and at most $16$ by \eqref{eq:local-diameter}. Since it equals
twice the diameter of $A_{\mathbb Z}$, it is $14$ or $16$.
Translate $A_{\mathbb Z}$ so that its minimum is zero. If the
diameter of $2A_{\mathbb Z}$ is $14$, then $2A_{\mathbb Z}=[0,14]$.
If it is $16$, then symmetry gives
$2A_{\mathbb Z}=[0,16]\setminus\{j,16-j\}$ with $1\le j<8$.
If $j\ge2$, the two omitted points equal $1+(j-1)$ and $1+(15-j)$,
respectively, with both summands in $2A_{\mathbb Z}$. The other
points of $[0,16]$ are sums with zero. Symmetry then gives
$[0,32]\subseteq4A_{\mathbb Z}$, contrary to $|4X|\le31$.
Thus $j=1$. The respective two-layer intervals are $[0,14]$ and
$[2,14]$, with losses $0$ and $4$. Their widths are at least the
diameter of $A_{\mathbb Z}$ and therefore of $U_{\mathbb Z}$.

\caseheading{4}{Layers of size five}
If there are fewer than two nonfirst layers of size five, leave
any such layer as a leftover. Otherwise take a balanced pair with
common support $U$, and set $X=\pi(\Sigma_2(U))$ and $D=X-X$.
Since $D=-D$ and $0\in D$,
its cardinality is odd. For $p\le31$, the local bound gives
$|D|\le13+e\le19$.

\caseheading{4.1}{$|D|\ge25$}
The pair provides the stated growth with $b=g=12$.

\caseheading{4.2}{$|D|\le23$}
Let $X^*=\Sigma_2(U)$ and $D^*=X^*-X^*$. The joint central sumset of the balanced pair is a translate of $D^*$, so Lemma~\ref{foundation:budget} gives injectivity on $D^*$ and hence on $X^*$, a translate of which lies in $D^*$. We show that $U$ admits integer coordinates with a unit common difference and that the difference set
$D_{\mathbb Z}=\Sigma_2(U_{\mathbb Z})-\Sigma_2(U_{\mathbb Z})$
satisfies the following bounds for some integer $a_0$.
\begin{equation}\label{eq:five-pure-goal}
 [-a_0,a_0]\subseteq D_{\mathbb Z}\subseteq[-h_0,h_0],\qquad
 a_0\ge6,\qquad 2a_0\ge\operatorname{diam}(U_{\mathbb Z}),\qquad
 h_0-a_0\le6,
\end{equation}
where $h_0=\max D_{\mathbb Z}$.

\caseheading{4.2.1}{$|X|\le8$}
The ten unordered pairs from $U$ give at most eight projected sums
in $C_p$. Since $\pi$ is injective on $X^*$, equal projected sums
are equal in $G$. Distinct pairs with equal sums are disjoint, and
no three pairs are pairwise disjoint. Thus there are at least two
pairing equalities in $G$. Each four-element subset gives at most
one such equality. Indeed, two different pairing equalities would
force two values to coincide, since $2$ is a unit in $G$. Choose
two such subsets sharing
three values $x,y,z$. Each fourth value is among $x+y-z$,
$x+z-y$, and $y+z-x$. Relabel the shared values so that the fourth
values are $x+y-z$ and $x+z-y$.
For $u=x-z$ and $v=y-z$ we obtain
$$
 U-z=\{0,v,u,u+v,u-v\},
 \qquad \overline u,\overline v\ne0.
$$

To count the formal sums before evaluation in $G$, use integer
coefficient pairs for $u,v$. The formal pair sums have
$v$-coefficients $\{1\}$, $[-1,2]$, and $[-1,1]$ in the rows
with $u$-coefficients $0,1,2$, respectively. Subtracting these
rows gives the following rows of the formal difference set
with nonnegative $u$-coefficient. The other rows follow by
negation.
$$
\begin{array}{c|ccc}
 u\text{-coefficient}&0&1&2\\ \hline
 v\text{-coefficients}&[-3,3]&[-3,2]&[-2,0]
\end{array}
$$
The formal difference set has $7+2(6+3)=25$ points. The
differences in any one row have distinct images under $\pi$.
Since $|D|\le23$, two points in different rows have the same image,
giving
$j\overline u=i\overline v$, where the possible $i$ are as follows.
$$
\begin{array}{c|cccc}
 j&1&2&3&4\\ \hline
 i&[-5,6]&[-4,6]&[-2,5]&[0,4]
\end{array}
$$
For example, the rows with coefficients $0,1$ give
$[-3,3]-[-3,2]=[-5,6]$, which contains the intervals from all
other pairs of adjacent rows. The other ranges follow by
the same subtraction.

Injectivity on $D^*$ gives $ju=iv$ in $G$. Let $b/a=i/j$ in
lowest terms with $a>0$. Since $a\le j\le4<p$, there is a unique
unit $d\in G$ with $ad=v$. Now $j(u-bd)=0$, so $u=bd$ since
$j<p$. Reflection by $t\mapsto v-t$ replaces $b/a$ by $1-b/a$,
so assume $b/a\ge1/2$. To justify using the same ranges after
reflection, note that before reflection we have $a\le4$ and
$|b|\le6$. The integer values of the formal differences have
absolute value at most $24$. As $48<p^2$, two such values are equal
exactly when the corresponding differences in $G$ are equal.
By injectivity on $D^*$, this is also equivalent to equality
after reduction modulo $p$. Reflection preserves the integer
difference set. After reflection two distinct rows still give
the same image, so the reflected ratio also comes from the four
ranges above.

Distinctness excludes $b/a\in\{-1,0,1,2\}$. The table gives the remaining ratios, with outer radii and omissions obtained from
$$
 D_{\mathbb Z}=a[-3,3]\ \cup\ \pm\bigl(b+a[-3,2]\bigr)\ \cup\ \pm\bigl(2b+a[-2,0]\bigr).
$$
Here $a[I]=\{ai:i\in I\}$.
$$
\begin{array}{c|r|c@{\qquad}c|r|c}
 (a,b)&h_0&[1,h_0]\setminus D_{\mathbb Z}&
 (a,b)&h_0&[1,h_0]\setminus D_{\mathbb Z}\\ \hline
 (1,3)&6 &\varnothing &(2,1)&6 &\varnothing\\
 (1,4)&8 &\varnothing &(2,3)&7 &\varnothing\\
 (1,5)&10&\varnothing &(2,5)&10&\varnothing\\
 (1,6)&12&\{9\}      &(3,2)&9 &\varnothing\\
 (3,4)&10&\varnothing &(3,5)&11&\varnothing\\
 (4,3)&12&\{10\}     &     &  &
\end{array}
$$
The inner radius is at least $6$ and $h_0-a_0\le4$. By \eqref{eq:support-difference}, $\operatorname{diam}(U_{\mathbb Z})\le h_0\le12$. This proves \eqref{eq:five-pure-goal}.

\caseheading{4.2.2}{$|X|\ge9$}
Here $p\ge29$, since for $p=23$ the bound in Lemma~\ref{foundation:budget}
gives $|D|\le15$, whereas Lemma~\ref{lem:cauchy-davenport} gives $|D|\ge17$. By Lemma~\ref{lem:local-integer-model} and \eqref{eq:planar-difference},
a rational model of dimension at least two would give
$|D|\ge3|X|-3\ge24$. Thus we may take $Y=\Sigma_2(U_{\mathbb Z})$
in the primitive integer model, where $|Y|=|X|$ and $|Y-Y|=|D|$.
Since $|D|-|X|\le14<2|X|-3$, \eqref{eq:local-diameter} gives
$h_0=\operatorname{diam}(Y)\le|D|-|X|\le14$.

We show that $[1,\max\{6,h_0-6\}]\subseteq Y-Y$. Translate so that $Y\subseteq[0,h_0]$ contains both endpoints. If $t>0$ is absent from $Y-Y$, then $Y$ contains no adjacent points in any residue class chain with step $t$. Write $h_0+1=jt+r$ with $0\le r<t$. There are $r$ chains of length $j+1$ and $t-r$ of length $j$, so
$$
 |Y|\le r\left\lceil\frac{j+1}{2}\right\rceil+(t-r)\left\lceil\frac j2\right\rceil.
$$
For $8\le h_0\le14$ and $1\le t\le\max\{6,h_0-6\}$, the right side can reach $9$ only in the following cases.
$$
\begin{array}{c|c|c}
 h_0&t&\text{Chain lengths}\\ \hline
 13&5&3,3,3,3,2\\
 14&3&5,5,5\\
 14&5&3,3,3,3,3\\
 14&6&3,3,3,2,2,2
\end{array}
$$
We rule out the four rows. Since $|Y-Y|\le23$, there are at most $11$ positive differences.

\begin{itemize}
\item For $h_0=13$ and $t=5$, $Y$ contains $[0,3]$, $[10,13]$, and one of $4,9$. Hence $Y-Y$ contains $[1,4]\cup[6,13]$.
\item For $h_0=14$ and $t=3$, $Y=[0,2]\cup[6,8]\cup[12,14]$. Hence $Y-Y$ contains $[1,2]\cup[4,8]\cup[10,14]$.
\item For $h_0=14$ and $t=5$, at least four chains contribute both endpoints. Thus $E\cup(E+10)\subseteq Y$ for a four-point set $E\subseteq[0,4]$. If $E$ omits an interior point, then $E-E=[-4,4]$ and $Y-Y$ contains $[1,4]\cup[6,14]$. Otherwise reflect to take $E=[0,3]$. Since $14\in Y$, we have $[10,14]\subseteq Y$, and the differences include $[1,4]\cup[7,14]$.
\item For $h_0=14$ and $t=6$, $Y$ contains $[0,2]$ and $[12,14]$, and one point from each of $\{3,9\},\{4,10\},\{5,11\}$. If it contains $9$ or $5$, the differences include $[1,5]\cup[7,14]$. Otherwise $3$ and $11$, together with the endpoint intervals, give $[1,3]\cup[8,14]$. The choice of $4$ or $10$ supplies $4$ and $7$, giving $[1,4]\cup[7,14]$.
\end{itemize}
Every row gives at least $12$ positive differences, a contradiction. Take $a_0=\max\{6,h_0-6\}$. Then $h_0-a_0\le6$ and $2a_0\ge h_0$. By \eqref{eq:support-difference}, $\operatorname{diam}(U_{\mathbb Z})\le h_0$.
This proves \eqref{eq:five-pure-goal}.

In both subcases of Case~4.2, the balanced pair gives an interval
of width $2a_0\ge12$, at least the support diameter. Its contribution
is $12$, and its loss is $2(h_0-a_0)\le12$. If $p\le31$, then
$|D|\le19$. For $|X|\le8$ the loss is already at most $8$. For
$|X|\ge9$ we have $h_0\le19-|X|\le10$, which again bounds the loss by $8$.

\par\medskip
\noindent\textbf{Aligned and leftover layers.}
If an initial interval is available, let $w$ be its width. The support
of each aligned layer is contained in the chosen support of size $q$
and has diameter at most $w$ in its integer coordinates. By
Lemma~\ref{lem:fixed-length-gaps}, the gaps between successive
central sums of such a layer are at most $w$. The same lemma shows
that adding any of these central sumsets to the initial interval
or to a wider interval again gives an interval. Each central sumset
has integer diameter at least its contribution $f(s_j)$.

A two-element layer is always aligned, since it provides an initial interval. The remaining unaligned layers are grouped as in the table. The weights of full groups of sizes three, four, and five are $32$, $28$, and $12$, respectively, each at most $16g/5$. The leftover weights sum to at most $3\cdot8+3\cdot7+6=51$.
\end{proof}

\section{Excluding the remaining multiplicities}\label{sec:residual}

\begin{theorem}\label{res:closure}
Let $p\ge11$. If $S$ is a sequence of $2p$ units, $\pi$ is injective on its support, and $m=p-\height(S)\ge m_0$, then $S$ has a subsequence $T$ with $|T|=p$ and $\sigma(T)\in pG\setminus\{0\}$.
\end{theorem}

\begin{proof}
Assume that $S$ is exceptional. Set $s=s_1$ and $u=s-4$, and let
$v_2,v_3$ be the second and third largest multiplicities.
Lemmas~\ref{foundation:budget} and~\ref{foundation:tail} give
\begin{equation}\label{res:tail-parameters}
 s\ge4,\qquad e+f(u)\le p-2m-3,\qquad
 v_2-e\ge m+1-e/2.
\end{equation}
For the last inequality, \eqref{eq:tail-weight-identity} gives
$$
 2(p+m)-2v_2\le2(p+m)-(v_2+v_3)\le b_{S}=2p-2-e.
$$
For $p=11,13$, the inequality $p-2m-3<0$ contradicts $e\ge0$.
Thus $p\ge17$. For $p=17,19$, \eqref{res:tail-parameters}
gives $e=0$, $s\le5$, and $v_2-e\ge m+1\ge8$.
For $p\ge23$, we have $4m\ge p+13$, so
$$
 2e+19\le2p-4m+13\le p,\qquad
 v_2-e\ge m+1-\frac{p-19}{4}\ge9.
$$

\medskip\noindent\emph{An integer interval from the layers other than the first.}

\caseheading{1}{$e=0$}
Lemma~\ref{foundation:tail} gives $\Phi\ge-p+9m+2>21$.
Thus some nonfirst layer has size two, or some size $r\in\{3,4,5\}$
occurs at least twice outside the first layer.

\caseheading{1.1}{Some nonfirst layer has size two}
Fix any central selection sizes of total $p$. The two-element layer
gives an interval of width $1$.

\caseheading{1.2}{No nonfirst layer has size two}
Choose two layers of size $r$ with common support $U$. Fix central
selection sizes of total $p$, taking $\lfloor r/2\rfloor$ and
$\lceil r/2\rceil$ in these two layers. Let
$B=\Sigma_{\lfloor r/2\rfloor}(U)$ and let
$Y=B+\Sigma_{\lceil r/2\rceil}(U)$. By \eqref{eq:sum-complement},
$\pi(Y)$ is a translate of $\pi(B)-\pi(B)$ when $r$ is odd and equals
$\pi(B)+\pi(B)$ when $r$ is even. Lemma~\ref{foundation:interpolation}
gives $|\pi(B)|\ge f(r)+1$.
By Lemma~\ref{foundation:budget} with $e=0$, we have
$|\pi(Y)|\le2f(r)+1\le p-2$. Lemma~\ref{lem:cauchy-davenport}
therefore gives
$$
 2f(r)+1\le2|\pi(B)|-1\le|\pi(Y)|\le2f(r)+1.
$$
All inequalities are equalities, so \eqref{eq:vosper} implies that
$\pi(B)$ is an arithmetic progression in $C_p$. Since $Y$ is $B+B$
or a translate of $B-B$, injectivity on $Y$ gives injectivity on
$B-B$ by the argument in Lemma~\ref{lem:small-layer-local}.
The successive differences of the progression therefore lift to
one common difference in $B$. Thus one layer provides an interval
of width $f(r)$.

In both subcases of Case~1, the interval comes from one layer and has width
equal to its contribution, at most $6$.

\caseheading{2}{$e>0$}
Then $p\ge23$ and Lemma~\ref{lem:small-layer-local} applies. Retain the
central selection sizes fixed in Section~\ref{sec:local-groups}.
For each full group indexed by $i$, let $g_i$ be its gain.
Suppose $\sum_i g_i>e$ and add groups until their total gain
first exceeds $e$. Let $J$ be the set of selected group indices,
$j_0$ the last one added, and $b$ their combined contribution.
For each $i\in J$, let $Y_i$ be the joint central sumset of the
$i$th group. Since each group contributes at most
its gain and at most $16$, we have
$$
 \sum_{i\in J\setminus\{j_0\}}g_i\le e<\sum_{i\in J}g_i,\qquad
 b\le\sum_{i\in J\setminus\{j_0\}}g_i+16\le e+16.
$$
Let $Y=\sum_{i\in J}Y_i$. Since $b+e+2\le2e+18\le p-1$
and $\sum_{i\in J}g_i\ge e+1$, Lemma~\ref{lem:cauchy-davenport} gives
$$
 |\pi(Y)|\ge\min\left\{p,1+b+\sum_{i\in J}g_i\right\}\ge b+e+2.
$$
The layers forming $Y$ have total contribution $b\le p-e-1$, so
Lemma~\ref{foundation:budget} gives $|\pi(Y)|\le b+e+1$, a
contradiction.
Thus $\sum_i g_i\le e$.

Let $\Phi_{\mathrm a}$ be the weight of the aligned layers. The full groups and leftovers have weight at most $16e/5+51$, so
$$
 \Phi_{\mathrm a}-4e>
 \frac{5p+107}{4}-\frac{11e}{5}-51
 \ge\frac{3p-67}{20}>0.
$$
Since the aligned layers have positive total weight, an initial
interval exists. For the initial interval and each remaining
aligned factor, let $w_i$ be its width and $b_i$ its contribution.
Each satisfies $w_i\ge b_i$, and each layer's weight is at most
four times its contribution. Thus
$$
 \sum_i(2w_i-b_i)\ge\sum_i b_i\ge\frac{\Phi_{\mathrm a}}4>e.
$$
Start with the initial interval and add aligned factors until
the sum on the left first exceeds $e$. By Lemma~\ref{lem:small-layer-local},
each factor's gaps are at most the initial width, so
Lemma~\ref{lem:fixed-length-gaps} ensures that the sum remains an interval.
The initial contribution is at most $12$. If further factors
are needed, the contribution before the final addition is at
most $e$ and the last factor contributes at most $6$. Thus the
total contribution is at most $\max\{12,e+6\}$.

In Cases~1 and~2, denote the resulting interval by $P$, its width by
$w$, its contribution by $b_P$, and its truncation loss by $l$.
Set $\gamma=w-b_P$. Then
$$
 \gamma\ge0,\qquad w+\gamma>e.
$$
The sum of the full integer diameters of the factors used, treating the
initial group as one factor, is $w+l$. Only the initial group is
truncated, and its retained width is at least its contribution. The
initial group uses at most four layers. When $e=0$, it uses one layer
and $l=0$.

\medskip\noindent\emph{Lifting the first layer.}
We first bound the combined contribution of the first and selected layers. For $p=17,19$, we have $f(s)+b_P\le12\le p-e-1$. For $p\ge23$, use $2u-f(u)\le4$ and \eqref{res:tail-parameters} to obtain
$$
 2u\le p-2m+1-e\le 2m-12-e.
$$
Since $2u$ is even, the last bound gives
$2u\le2m-12-2\lceil e/2\rceil$. The bounds $b_P\le6$ for $e=0$
and $b_P\le\max\{12,e+6\}$ for $e>0$ imply
$b_P+2\le12+2\lceil e/2\rceil$. Thus $2u+b_P+2\le2m$, and
\eqref{res:tail-parameters} gives
$$
 f(s)+b_P=f(u)+(2u+b_P+2)+2\le p-e-1.
$$

Let $A$ be the central sumset of the first layer in $G$, and
regard $P$ here as the set $\{y+id:i\in[0,w]\}\subseteq G$,
where $y$ is the element of $P$ at coordinate zero and $d$ is
its unit common difference. By
Lemma~\ref{foundation:budget}, $\pi$ is injective on
$A+P$ and $|\pi(A+P)|\le f(s)+b_P+e+1$. Since a translate of $A$
lies in $A+P$, $\pi$ is injective on $A$. List $a_i\in A$ in the
cyclic order of their images along $\overline d$. Let $\ell_i$ be
the positive distance from $\pi(a_i)$ to $\pi(a_{i+1})$ in this
order, so $\sum_i \ell_i=p$. The coverage in each gap is
$\min\{\ell_i,w+1\}$, hence
$$
 |\pi(A+P)|=\sum_i\min\{\ell_i,w+1\}
           =|\pi(A)|+\sum_i\min\{\ell_i-1,w\}.
$$
Together with $|\pi(A)|\ge f(s)+1$, this gives
$$
 \sum_i\min\{\ell_i-1,w\}\le b_P+e.
$$
Whenever $\ell_i\le w$, the elements $a_i+y+\ell_i d$ and $a_{i+1}+y$
of $A+P$ have the same image. Injectivity on $A+P$ gives
$a_{i+1}-a_i=\ell_i d$. If every $\ell_i\le w$, summing these equalities
cyclically gives $0=(\sum_i \ell_i)d=pd$, a contradiction.

Thus some gap $\ell_{i_0}$ is at least $w+1$. It contributes $w$ to
the last sum, and the other contributions sum to at most
$b_P+e-w=e-\gamma<w$. In particular, each other gap is at most $w$.
Cutting the cyclic order of $\pi(A)$ at $\ell_{i_0}$ gives integer
coordinates along $\overline d$ with diameter
$$
 \sum_{i\ne i_0}\ell_i=|\pi(A+P)|-w-1\le f(s)+e-\gamma=:L<p.
$$
The overlaps at the other gaps give the same successive
differences in $A$, so $A$ itself lifts to that diameter.
Since $|\pi(A)|\ge2$, the large gap is at most $p-1$,
and $w\le p-2$.

By \eqref{eq:support-difference}, $L_1-L_1\subseteq A-A$. Choose
$c\in L_1$ and, for each $b\in L_1$, choose $x_b\in[-L,L]$ with
$b-c=x_bd$ in $G$. For any $b,b'\in L_1$, we have
$$
 x_b-x_{b'}\equiv t\pmod{p^2},\qquad t\in[-L,L].
$$
Since $|x_b-x_{b'}-t|\le3L<p^2$, we have $x_b-x_{b'}=t$.
Thus the coordinates have diameter at most $L$. Let $a$ be the
support value of least coordinate. After translating the
coordinates so that $a$ has coordinate zero, the support is
$$
 \{a+x_i d:1\le i\le s\},\qquad 0=x_1<\cdots<x_s=D\le L.
$$
In these coordinates, a selection of $p$ terms with integer sum
$z$ has sum $pa+zd$ in $G$.
The difference between two integer central sums has absolute value
at most $sL$. Its image after multiplication by $d$ lies in
$A-A$, so the difference is congruent modulo $p^2$ to an integer
in $[-L,L]$.
Since $s\le p-1$ and $L<p$, the bound $(s+1)L<p^2$ makes the
two integers equal. The central sumset of the first layer
therefore has integer diameter at most $L$.

For a layer of size $r$ with coordinates $y_1<\cdots<y_r$, the width of its central sumset is
\begin{equation}\label{res:width-gap}
 f(r)+\sum_{i=1}^{r-1}\min\{i,r-i\}(y_{i+1}-y_i-1).
\end{equation}
To see \eqref{res:width-gap}, note that the two central selection
sizes give the same width by \eqref{eq:sum-complement}. Set $k=\lfloor r/2\rfloor$.
Subtract the minimum central sum from the maximum and
collect terms by gaps to obtain
$$
 \sum_{i=1}^{k}(y_{r-k+i}-y_i)
   =\sum_{i=1}^{r-1}\min\{i,r-i\}(y_{i+1}-y_i).
$$
The sum of the gap weights is $k(r-k)=f(r)$, which gives
\eqref{res:width-gap}. For the first layer the total excess
above $f(s)$ is at most $L-f(s)=e-\gamma$. Every gap
weight is at least one, so
$$
 D\le s-1+e-\gamma<p-1,\qquad
 \max_i(x_{i+1}-x_i)\le e-\gamma+1\le w.
$$
Here the strict diameter bound follows from
$s+e\le p-2m+1+u-f(u)\le p-2m+2$.

\medskip\noindent\emph{Completion to $p$ terms.}
Let $\eta_j$ be the excess of the width in \eqref{res:width-gap}
over $f(s_j)$, and set $\eta=\sum_j\eta_j$. Let $k$ be the number of
terms supplied by the selected layers. The fixed central selection
sizes sum to $p$, so the central selections in the remaining layers
supply $p-k$ terms. Let $z,z'$ be the sums of a selection from
the selected layers in the coordinates used to construct $P$
and in the full support coordinates, respectively. Let $z_0,z'_0$
be the sums of a fixed selection in the same two coordinate
systems. The selection sizes in each layer are fixed, so the
origins cancel when subtracting the two group sums in either
coordinate system. Thus $(z-z_0)d=(z'-z'_0)d$ in $G$. Since $d$
has order $p^2$ and since $|z-z_0|\le w+l$,
$|z'-z'_0|\le kD\le pD$, $D\le p-2$, $w+l\le p+10$, and
$p\ge17$, we obtain
$$
\begin{aligned}
 (z-z_0)-(z'-z'_0)&\equiv0\pmod{p^2},\\
 |(z-z_0)-(z'-z'_0)|&\le pD+w+l\le p(p-2)+p+10<p^2.
\end{aligned}
$$
It follows that $z'=z+(z'_0-z_0)$ in the integers. The same
translation takes the full sumset and the coordinate interval
of $P$ into the full support coordinates, preserving their
widths $w+l$ and $w$ and hence the loss $l$. Denote the resulting
integer interval by $P_{\mathbb Z}$.

All first-layer positions remain outside the selected layers, so the
remaining positions contain every support value and have maximum gap
at most $w$. By Lemma~\ref{lem:fixed-length-gaps}, adding their
integer sumset of $p-k$ terms to $P_{\mathbb Z}$ gives an interval. The central selections
in the remaining layers show that its width is at least
$b_{S}+\eta-l$. If $\eta>e+l$, then
$$
 b_{S}+\eta-l\ge b_{S}+e+1=2p-1.
$$
The interval then contains $jp$ and $(j+1)p$ for some integer $j$.
The corresponding sums in $G$ are $pa+jpd$ and $pa+(j+1)pd$.
They belong to $\Sigma_p(S)\cap pG$ and differ by $pd\ne0$,
contrary to \eqref{eq:exceptional}. Thus $\eta\le e+l$.

The coordinates of the two most frequent values lie in the same
maximal interval $[v_-,v_+]$ of consecutive support coordinates.
Otherwise each of the first $v_2$ layers has a gap
and contributes at least one unit of excess, except for the layers
in the initial group, which we count separately. Fixing the selections
in the other selected layers shows that the coordinate comparison
preserves the width of the full central sumset of the initial group.
Its full width exceeds its retained width by $l$, and its retained
width is at least its contribution. Thus the initial group has total
excess at least $l$. Since it uses at most four layers, we get
$\eta-l\ge v_2-4>e$, a contradiction.

Let $b_{\mathrm{out}}$ be the number of terms with coordinates
outside $[v_-,v_+]$, and let $q_j$ be the number of such values in
layer $j$. Set
$$
 E=D-(s-1),\qquad M_0=\sum_j\min\{2,q_j\}.
$$
The layers are nested. A nonempty layer with $j>v_2$ consists
only of the most frequent value, whose coordinate lies in
$[v_-,v_+]$. Thus $q_j>0$ implies $j\le v_2$, and the layer
contains the two most frequent values. In layer $j$, let $q_-$
and $q_+$ be the counts to the left and right of $[v_-,v_+]$,
respectively. Thus $q_-+q_+=q_j$. If a side contains $q>0$
coordinates, its boundary gap has size at least two and weight
$\min\{q,s_j-q\}\ge\min\{2,q\}$ in \eqref{res:width-gap}.
The two sides use distinct gaps when both are present, so
$$
 \eta_j\ge\min\{2,q_-\}+\min\{2,q_+\}\ge\min\{2,q_j\}.
$$
For the first layer the total gap excess is $E$. After
allocating one unit of gap at each boundary to the preceding
estimate, at least $E-2$ units remain. Therefore
$\eta_1\ge\min\{2,q_1\}+E-2$ and $E+M_0\le\eta+2$.

For $q_j\ge3$ we have
$$
 2(q_j-2)\le f(q_j)\le f(s_j-2)=f(s_j)-s_j+1.
$$
Since $f(s_j)\ge s_j-1$ for every nonempty layer, summing yields
$$
 \begin{aligned}
 b_{\mathrm{out}}
 &=M_0+\sum_j\max\{q_j-2,0\}\\
 &\le M_0+\tfrac12\bigl(b_{S}-p-m\bigr),
 \end{aligned}
$$
where $\sum_j(s_j-1)=p+m$. Together with $\eta\le e+l$
and $b_{S}=2p-2-e$, this gives
\begin{equation}\label{res:span-outliers}
 D+b_{\mathrm{out}}\le s+\frac{p-m+e}{2}+l
 \le p-\frac{3m-9-2l}{2}<p.
\end{equation}
The second inequality uses \eqref{res:tail-parameters} and $2u-f(u)\le4$. For strictness when $e=0$, use $m\ge7$ and $l=0$. If $e>0$ and $p\le31$, then $m\ge9$ and $l\le8$. If $p\ge37$, then $m\ge13$ and $l\le12$.

Fix any subsequence $T$ whose terms have coordinates outside
$[v_-,v_+]$ and complete it with $q=p-|T|$ terms whose coordinates
lie in that interval. Let $z_1\le\cdots\le z_n$ be the ordered
coordinates of all terms with coordinates in $[v_-,v_+]$, with
$n=2p-b_{\mathrm{out}}$. We have
$p-b_{\mathrm{out}}\le q\le p\le n$ by \eqref{res:span-outliers}.
Since $z_{i+1}-z_i\le1$ for $1\le i<n$,
Lemma~\ref{lem:fixed-length-gaps} shows that the integer sums of
the completed selections form an interval $I_T$. For
$k=\min\{q,n-q\}$, its width is
$$
 \operatorname{diam}(I_T)=\sum_{i=1}^k(z_{n-k+i}-z_i).
$$
A summand can vanish only at a value with multiplicity greater
than $n-k\ge n/2$. At most one value has that multiplicity,
and the number of zero summands is at most
$\max\{0,\height(S)-(n-k)\}$. We have
$k\ge p-b_{\mathrm{out}}$ from $q\ge p-b_{\mathrm{out}}$ and
$n-q\ge n-p=p-b_{\mathrm{out}}$. Hence
$$
 \operatorname{diam}(I_T)\ge
 \min\{k,n-\height(S)\}\ge p-b_{\mathrm{out}}>D.
$$
If $T'$ is obtained by adding one term with coordinate outside
$[v_-,v_+]$ to $T$, then at each
endpoint one coordinate from $[v_-,v_+]$ is replaced by the
coordinate of the added term. All coordinates lie in $[0,D]$, so
$$
 |\min I_{T'}-\min I_T|\le D,\qquad
 |\max I_{T'}-\max I_T|\le D.
$$
Since both intervals have width greater than $D$, they overlap.
All choices of $T$ are connected by successive additions or
removals, so their union is an interval $I$ of all $p$-term
integer sums.

The interval $I$ contains all central sums, so
$$
 \operatorname{diam}(I)\ge b_{S}\ge\sum_j(s_j-1)=p+m\ge p-1.
$$
Choose $pj\in I$. If $a+jd$ is a unit, the corresponding
$p$-term subsequence has sum $p(a+jd)\in pG\setminus\{0\}$.
Otherwise $j\notin[v_-,v_+]$, since every coordinate in this
interval represents a unit. Let $t\in I$ be an integer sum of
$p$ terms with coordinates in $[v_-,v_+]$.
Then
$$
 pj<pv_-\le t\qquad\text{or}\qquad t\le pv_+<pj.
$$
Thus $pv_-\in I$ in the first case and $pv_+\in I$ in the
second. The corresponding endpoint represents a unit,
again giving a $p$-term sum in $pG\setminus\{0\}$.
\end{proof}

\section{Proofs of the main results}

\begin{proof}[Proof of Theorem~\ref{thm:main}]
Assume first that $S$ is exceptional. Since any subsequence $a^p$
would have sum $pa\in pG\setminus\{0\}$, we have $\height(S)\le p-1$.

\caseheading{1}{$\pi$ is not injective on $\supp(S)$}
Proposition~\ref{prop:collision} gives the first form.

\caseheading{2}{$\pi$ is injective on $\supp(S)$}
Let $m=p-\height(S)$ and let $v_2,v_3$ be the second
and third largest multiplicities. Since $v_2+v_3\le2(p-m)$,
Lemma~\ref{foundation:budget} and \eqref{eq:tail-weight-identity} give
$$
 4m\le2(p+m)-(v_2+v_3)\le b_{S}\le2p-2.
$$
Thus $1\le m\le(p-1)/2$. If $p=3$, injectivity on $\supp(S)$
gives $|\supp(S)|\le2$, so $|S|\le4<2p$. Hence $p\ge5$.

\caseheading{2.1}{$m=1$}
By Proposition~\ref{prop:base}, there are two distinct values $a,b$
of multiplicity $p-1$. Hence $S=a^{p-1}b^{p-1}xy$ with
$x,y\notin\{a,b\}$, and Proposition~\ref{prop:double} gives the
second form.

\caseheading{2.2}{$m\in\{2,3\}$}
Proposition~\ref{prop:23} gives a contradiction.

\caseheading{2.3}{$m\ge4$}
Since $m\le(p-1)/2$ and $p$ is prime, we have $p\ge11$.
The definition of $m_0$ gives
$$
 m<m_0
 \quad\Longleftrightarrow\quad
 m\le\left\lfloor\frac p3\right\rfloor+1
 \quad\text{and}\quad
 m\le\left\lfloor\frac{p+11}{4}\right\rfloor
 \quad\Longleftrightarrow\quad
 p\ge\max\{3m-3,4m-11\}.
$$
If $m<m_0$, Proposition~\ref{thm:high} gives a contradiction.
If $m\ge m_0$, Theorem~\ref{res:closure} gives a contradiction.

Conversely, the calculation in the proof of
Proposition~\ref{prop:collision} shows that for the first form the
only $p$-term subsequences with sum in $pG$ are $A$ and $B$.
Both have sum zero. For the second form, reducing the formulas for
$x$ and $y$ modulo $p$ gives
$$
 \begin{aligned}
 \overline x&=\overline b-i(\overline a-\overline b),\\
 \overline y&=\overline b-(p-1-i)(\overline a-\overline b).
 \end{aligned}
$$
Since $\overline a\ne\overline b$ and $1\le i,p-1-i\le p-2$,
both $\overline x$ and $\overline y$ lie outside
$\{\overline a,\overline b\}$. Moreover,
$\overline x-\overline y=(p-1-2i)(\overline a-\overline b)$.
Since $|p-1-2i|\le p-3$, equality of $\overline x$ and $\overline y$
forces $2i+1=p$, in which case $x=y$.
Thus $\pi$ is injective on the support, and
Proposition~\ref{prop:double} proves sufficiency.
\end{proof}

\begin{proof}[Proof of Corollary~\ref{cor:direct}]
Suppose $W$ is a sequence of $2p+1$ units with no $p$-term
subsequence whose sum belongs to $pG\setminus\{0\}$. For each term $w$
of $W$, the sequence $Ww^{-1}$ is exceptional. By Theorem~\ref{thm:main},
its largest residue class contains $p$ terms in the first form and
$p-1$ in the second. Indeed, in the second form the classes of $a$ and
$b$ each contain $p-1$ terms, and any other class contains at most
two terms. If some class of $W$ contained at least $p+1$ terms, we
could retain $p+1$ such positions by deleting another position.
Thus every class of $W$ contains at most $p$ terms.

\caseheading{1}{Some residue class contains exactly $p$ terms}
Let $A$ consist of the $p$ terms in this class and set $B=WA^{-1}$.
For every term $b$ of $B$, the sequence $Wb^{-1}=A(Bb^{-1})$
retains $A$ and therefore has the first form. Thus the terms of
$Bb^{-1}$ all lie in one residue class. Any two positions of $B$
survive when a third position is deleted, so all $p+1$ terms of
$B$ lie in a single residue class, a contradiction.

\caseheading{2}{Each residue class contains at most $p-1$ terms}
If $p=3$, the two nonzero classes contain at most four terms, a
contradiction. Hence $p\ge5$. Every deletion has the second form
and gives two classes with $p-1$ terms each. Since
$3(p-1)>2p+1$, exactly two classes in $W$ contain $p-1$ terms.
Deleting a term from either of them leaves at most one such
class, a contradiction.

Both cases give a contradiction. Therefore $\mathsf{s}_p^*\le2p+1$.

For the reverse inequality, let $A=1^{p-1}(1-p)$ and
$B=(-1)^{p-1}(p-1)$. Then $|A|=|B|=p$,
$\overline A=1^p$, $\overline B=(-1)^p$, and
$\sigma(A)=\sigma(B)=0$. Thus $AB$ has the first form and is
exceptional, which gives $\mathsf{s}_p^*\ge2p+1$.
\end{proof}

\section*{Declaration of generative AI and AI-assisted technologies}
The authors used a large language model during the preparation of this manuscript
for assistance with mathematical exploration, checking arguments,
organization and language editing. The authors are responsible for
the mathematical statements, proofs, references and final text.

\bigskip
\noindent
\begin{minipage}{\linewidth}
\small
\hypertarget{author-cong-lin}{Cong Lin}\\
Center for Combinatorics,
Nankai University,
Tianjin, P. R. China\\
E-mail: \href{mailto:9820260111@nankai.edu.cn}{\texttt{9820260111@nankai.edu.cn}}
\end{minipage}

\end{document}